\documentclass[12pt,reqno]{amsart}

\usepackage[dvips]{graphicx} 

\usepackage{
hyperref
}

\usepackage{breakurl}   

\usepackage{epigraph}

\theoremstyle{definition}

\numberwithin{equation}{section}

\newcommand\N {{\mathbb N}} 

\newcommand\R {{\mathbb R}}

\newcommand\astr{{{}^\ast\hspace{-1pt}\R}}

\title [Historical infinitesimalists and modern historiography]
       {Historical infinitesimalists and modern historiography of
         infinitesimals}

\author[J. B.]{Jacques Bair}\address{HEC-ULG, University of Liege,
  4000 Belgium}\email{j.bair@ulg.ac.be}

\author[A. B.]{Alexandre Borovik} \address{School of Mathematics,
  University of Manchester, Oxford Street, Manchester, M13 9PL, United
  Kingdom} \email{alexandre@borovik.net}

\author[V. K.]{Vladimir Kanovei} \address{IITP RAS, Moscow,
  Russia}\email{kanovei@googlemail.com}

\author[M. K.]{Mikhail G. Katz}\address{Department of Mathematics, Bar
  Ilan University, Ramat Gan 5290002
  Israel}\email{katzmik@math.biu.ac.il}

\author[S. K.]{Semen S. Kutateladze}\address{Sobolev Institute of
  Mathematics, Novosibirsk State University, Russia}
\email{sskut@math.nsc.ru}

\author[S. S.]{Sam Sanders} \address{Department of Philosophy 2, RUB
  Bochum, Bochum, Germany
  \url{http://sasander.wix.com/academic}}\email{sasander@me.com}

\author[D. S.]{David Sherry} \address{Department of Philosophy,
  Northern Arizona University, Flagstaff, AZ 86011, US}
\email{David.Sherry@nau.edu}

\author[M. U.]{Monica Ugaglia} \address{Il Gallo Silvestre, Localit\`a
  Collina 38, Montecassiano, Italy}\email{monica.ugaglia@gmail.com}


\begin{document}

\begin{abstract}
In the history of infinitesimal calculus, we trace innovation from
Leibniz to Cauchy and reaction from Berkeley to Mansion and beyond.
We explore 19th century infinitesimal lores, including the approaches
of Sim\'eon-Denis Poisson, Gaspard-Gustave de Coriolis, and
Jean-Nicolas No\"el.  We examine contrasting historiographic
approaches to such lores, in the work of Laugwitz, Schubring, Spalt,
and others, and address a recent critique by Archibald et al.  We
argue that the element of contingency in this history is more
prominent than many modern historians seem willing to acknowledge.
\end{abstract}


\thispagestyle{empty}



\maketitle
\tableofcontents


\section{Debate over Leibniz}
\label{s1}

Both Leibniz and Cauchy used the term \emph{infinitesimal} in their
work.  The meaning of the term has been the subject of scholarly
debates.  In this section, we focus on Leibniz's use of the term.  In
Sections~\ref{s2} and~\ref{s3}, we focus on Cauchy and his
interpreters.  In Section~\ref{s20}, we analyze certain
historiographic assumptions underlying existing interpretations of
these pioneers of infinitesimal analysis.  In Section~\ref{s1v}, we
address a recent critique by Archibald et al.  We summarize our
conclusions in Section~\ref{s6}.

George Berkeley claimed to find shortcomings in both the Newtonian and
the Leibnizian calculus.  While modern scholars (both historians and
mathematicians) also find shortcomings, their quibbles about Leibniz
are not identical to Berkeley's.  Berkeley's empiricism obscured from
him the coherence of the procedures of the Leibnizian infinitesimal
calculus.  Specifically, Berkeley's logical criticism overlooked the
coherence of Leibniz's relation of infinite proximity,
%
%
as we detail in Section~\ref{s81}.

\subsection{Berkeley's criticisms; Transcendental Law of Homogeneity}
\label{s81}

Berkeley was an English cleric whose \emph{empiricist} (i.e., based on
sensations) metaphysics tolerated no conceptual innovations, like
infinitesimals, without an empirical counterpart.  Berkeley was
similarly opposed, on metaphysical grounds, to infinite divisibility
of the continuum (which he referred to as \emph{extension}), an idea
widely taken for granted today (as it was already by Leibniz).

In addition to his \emph{metaphysical criticism} of the infinitesimal
calculus of Newton and Leibniz, Berkeley put forth a \emph{logical
  criticism} in his pamphlet \emph{The Analyst}.  He claimed to have
detected a logical fallacy at the foundation of the method.  The
distinction between logical and metaphysical criticisms in Berkeley
goes back to Sherry's 1987 article \cite{Sh87}; see further in
\cite{13a}.  In terms of Fermat's technique of adequality%
\footnote{Such a technique was used by Fermat to solve problems of
  finding tangents to curves, maxima and minima, and others; see
  further in \cite{13e} and \cite{18d}.}
exploiting an increment~$E$, Berkeley's objection can be formulated as
follows: the increment~$E$ is assumed to be nonzero at the beginning
of the calculation, but zero at its conclusion, an apparent logical
fallacy.

However, as noted by Fermat historian Str{\o}mholm
\cite[p.\;51]{St68},~$E$ is not assumed to be zero at the end of the
calculation, but rather is \emph{discarded} at the end of the
calculation.  Such a procedure was the foundation of both Fermat's
adequality and Leibniz's Transcendental Law of Homogeneity (TLH),
involving the relation of infinite proximity.  Leibniz discussed the
TLH in texts from 1695%
\footnote{See Bella \cite[p.\,192, note~70]{Be19}.}
and in a 1710 text \cite{Le10b}.%
\footnote{The 1710 text was analyzed by Bos \cite{Bo74} (see p.\;33
  and note~64 there); see further in~\cite{12e}.  See Bascelli et
  al.~\cite{Ba14} and Dinis \cite{Di22} for a study of modern
  formalisations.}
The TLH is closely related to a pair of modern procedures%
\footnote{On the dichotomy of procedures \emph{vs} foundations, see
  Section~\ref{s97}.}
in analysis:
\begin{enumerate}
\item
passing to the \emph{limit} of a typical expression such as
$\frac{f(A+E)-f(A)}{E}$ in the Weierstrassian approach, and
\item
taking the \emph{standard part} in infinitesimal analysis \cite{Ro66}.
\end{enumerate}

Meanwhile, Berkeley's own attempt to explain the calculation of the
slope when 
$y=x^2$ in Section XXIV of \emph{The\;Analyst} contains a logical
circularity.  Namely, Berkeley's argument relies on the determination
of the tangents of a parabola by Apollonius (which is equivalent to
the calculation of the slope).  The circularity in Berkeley's argument
is analyzed by Andersen \cite{An11}.  Far from exposing logical flaws
in the Leibnizian calculus, Berkeley's \emph{The\;Analyst} is itself
logically flawed.  Berkeley's character has been analyzed by Moriarty
\cite{Mo21} and \cite{Mo22}.

Berkeley's rhetorical flourishes such as the \emph{ghosts of departed
  quantities} were popular with a generation of scholars who
attributed exaggerated significance to his influence in the history of
the calculus.  These are the historians and Leibniz scholars of the
period until around 1966, including Boyer and Kline.  These scholars
\begin{enumerate}
\item
believed Berkeley to have provided the motivation for the eventual
success of the ``great triumvirate'' \cite[p.\;298]{Bo49} of Cantor,
Dedekind, and Weierstrass in eliminating the ghosts that haunted the
early calculus of Newton and Leibniz; and
\item
sought to
interpret Cauchy as anticipating the Weierstrassian \emph{Epsilontik}
with its alternating quantifiers.%
\footnote{The \emph{modus operandi} of such scholars can therefore be
  described in terms of a quest for the \emph{ghosts of departed
    quantifiers}; see \cite{17a}.}
\end{enumerate}

In short, Weierstrass-trained historians tended to attribute special
significance to Berkeley's critique because the Weierstrassian real
line is taken to exhibit purely Archimedean behavior admitting no
relation of infinite proximity.  But such an approach to the
historical infinitesimal calculus risks walking royal roads; see
Section~\ref{s20}.  Received attitudes underwent a subtle but
significant change in the second half of the 20th century.

\subsection{Changing attitudes toward Berkeley}
\label{s12}

Abraham Robinson developed modern infinitesimal analysis in his 1966
book \cite{Ro66}, building upon earlier work by Skolem \cite{Sk33},
Hewitt \cite{He48}, {\L}o\'s \cite{Lo55}, and others.  Robinson named
his theory ``Non-standard Analysis since it involves and was, in part,
inspired by the so-called Non-standard models of Arithmetic whose
existence was first pointed out by T. Skolem'' \cite[p.\,vii]{Ro66}.
Attitudes among scholars toward Berkeley's criticisms have undergone a
perceptible change since the appearance of Robinson's book.%
\footnote{Bockstaele's article published in the same year (1966)
  presents the two sides of the Belgian debate over the teaching of
  infinitesimals (see Section~\ref{s24} below),
%
%
but chooses to describe only one side as ``obstinate'' \cite[pp.\;2,
  8]{Bo66}.  No\"el, an advocate of infinitesimals, is described
mockingly as a ``never-desponding defender of the infinitely small''
\cite[p.\,14]{Bo66}.  Bockstaele concludes by mentioning ``a
definitive acknowledgment of the limit concept as the foundation of
the calculus'' \cite[p.\,16]{Bo66}.  See Section~\ref{s20} for an
analysis of teleological aspects of this type of historical
scholarship.}
The current generation of Leibniz scholars holding received opinions
derives its impetus from the second, 1990 edition of Ishiguro's book
\cite{Is90} on Leibniz, where she wrote:
\begin{enumerate}\item[]
Robinson's success in introducing infinitesimals into the
Weierstrassian analysis seemed to vindicate Leibniz from Berkeley's
famous attack, in which Berkeley claimed ``to conceive a quantity
infinitely small, that is, infinitely less than any sensible or
imaginable quantity or any finite quantity however small is, I
confess, beyond my capacity.''  \cite[p.\;83]{Is90}
\end{enumerate}
Ishiguro went on to disagree with such a seeming ``vindication,''
arguing that Leibniz was more rigorous than historians believed him to
be, in the following sense.  Ishiguro and her followers claim that
Leibniz never thought of infinitesimals as mathematical entities in
the first place, and that occurrences of the term \emph{infinitesimal}
in Leibniz do not \emph{refer} to a mathematical entity; they are mere
stenography for a more long-winded argument \`a la Archimedean
exhaustion.%
\footnote{Some Leibniz scholars are beginning to re-evaluate the
  claims of the Ishiguro school; see e.g., Esquisabel and Raffo
  Quintana \cite{Es21} and Samuel Henry Eklund at the University of
  California at Irvine \cite{Ek20}.  See further in \cite{21a} and
  \cite{23d}.}
Such scholars are mostly silent as to the pertinence of Berkeley's
critique.  What about Berkeley the great influencer?  One finds little
about this other than in reprints of the old classics by Boyer
\cite{Bo49} and Kline \cite{Kl80}.%
\footnote{Boyer and Kline also jointly fabricated the ``quote'' from
  Cavalieri about rigor being the concern of philosophy rather than
  mathematics.
}

\subsection{Potential infinity \emph{vs} infinite wholes; \emph{infinita
    terminata}}
\label{s13}

An aspect of Leibniz's thought at variance with modern usage is his
rejection of infinite wholes.  The dichotomy of potential infinity
\emph{vs} infinite wholes can be traced back to Aristotle, for whom
potential infinity meant an iterative process of repeating a procedure
over and over again ($\alpha\varepsilon\iota$).  Thus, one can only
have finite lines, but one can envision a process of doubling them
each time, again and again (similarly, dividing a line segment again
and again will give a smaller and smaller segment).  This is the
meaning of the so-called \emph{syncategorematic} infinite.  But the
process never leads to an infinite whole (for details see
Ugaglia~\cite{Ug22}).

Leibniz agreed with this conclusion.  His analysis of the Galilean
paradox (comparing integers and squares of integers) led him to reject
infinite wholes as contradictory,%
\footnote{\label{f8}Cauchy reached similar conclusions; see main text
  at note~\ref{f19}.}
and more precisely contradicting Euclid's part-whole principle
(whether or not Euclid himself meant for the principle to apply to
infinite pluralities is a separate issue); for details see
\cite[Section~4.1]{23d}, \cite{22b}.

In this connection, Leibniz elaborated an important distinction in
Proposition\;11 of his \emph{De Quadra\-tura Arithmetica}
\cite[pp.\;520--676]{Le2012}.%
\footnote{A French translation by Parmentier is available
  \cite{Le04b}.}
This is the distinction between bounded infinity (\emph{infinitum
  terminatum}) and unbounded infinity (\emph{infinitum interminatum}).
The latter, exemplified by an unbounded infinite line, is a
contradictory notion.  The former, exemplified by a segment with
infinitely separated endpoints, is a useful concept in geometry and
calculus; see further in \cite[Section~2.2]{21a}.  The distinction
between bounded and unbounded infinity was mentioned in the
29\;july\;1698 letter to Bernoulli \cite[III 523]{Ge50} as well as the
2 february 1702 letter to Varignon \cite[p.\,91]{Le02}.  The same
letter to Varignon contains a definition of an infinitesimal as a
``fraction infiniment petite, ou dont le denominateur soit un nombre
infini''%
\footnote{``An infinitely small fraction, or one whose denominator is
  an infinite number.''}
\cite[p.\;93]{Le02}, i.e., the reciprocal of an \emph{infinitum
  terminatum}.

It is therefore difficult to agree with Gert Schubring's claims that
\begin{enumerate}
\item
``{}`the founders of the calculus' had not at all created it with a
  non-Archime\-dean continuum in mind'' \cite[p.\;6]{Sc22}, or that
\item
``Leibniz always refused to be identified with a foundation based
  on---rather vaguely conceived---infinitesimals'' \cite[p.\,1]{Sc22}.
\end{enumerate}
While Schubring apparently believes infinitesimals to be ``vaguely
conceived,'' Leibniz defined them as fictional inassignable quantities
smaller than any assignable quantity, or magnitudes incomparable
with~$1$ in the sense of the violation of Definition 4 of Euclid's
Book V; see further in~\cite{21a} and \cite{18a}.%
\footnote{\label{f3}Such approaches to infinitesimals admit
  straightforward formalisations in modern infinitesimal theories: an
  infinitesimal is a number smaller in absolute value than every
  positive standard number; for details see e.g., \cite{17f} or
  \cite{21e} or \cite{23b}.  The viability of the application of
  nonstandard analysis to interpreting the procedures of the
  historical infinitesimalists depends on the procedure/foundation
  distinction; see Section~\ref{s97}.  On Schubring see also
  note~\ref{f7}.}

\subsection{Magnitudes \emph{vs} multitudes; Maxima and Minima}
\label{s14}

The significance of the Leibnizian distinction between bounded and
unbounded infinity tends to be underplayed by Leibniz scholars who
follow Ishiguro~\cite{Is90} in seeking to interpret Leibnizian
infinitesimals as stenography for exhaustion procedures in the sense
outlined in Section~\ref{s13}, and to relate such an approach to the
Scholastic concept of syncategorematic infinity (closely related to
potential infinity).  But while syncategorematic infinity is a
well-known and important concept in Leibnizian thought, the locution
\emph{syncategorematic infinitesimal} is nowhere to be found in known
Leibnizian texts.  Leibniz's work allows for an interpretation that he
worked with infinitesimals and bounded infinities as fictional
mathematical entities, rejecting infinite wholes while adhering to the
part-whole principle; see further in \cite{21a}.

Leibniz used the term \emph{Maxima} to refer to infinite wholes, and
the term \emph{Minima} to refer to points viewed as constituent parts
of the continuum.  Leibniz rejected both Maxima and Minima in the
following terms:
\begin{enumerate}\item[]
\emph{Scholium}.  We therefore hold that two things are excluded from
the realm of intelligibles: minimum and maximum; the indivisible, or
what is entirely \emph{one}, and \emph{everything}; what lacks parts,
and what cannot be part of another.%
\footnote{``Habemus ergo exclusa rebus intelligibilibus, duo: Minimum
  et Maximum; indivisibile, vel omnino unum, et omne; quod partibus
  careat, et quod pars alterius esse non possit'' (A VI 3 98).}
(Leibniz as translated by Arthur in \cite[p.\,13]{Ar01})
\end{enumerate}
Leibniz's rejection of Maxima amounts to the rejection of infinite
wholes (e.g., unbounded lines) as inconsistent.  The rejection of
their counterparts, Minima, amounts to the rejection of putative
simplest constituents of the continuum, i.e., the rejection of a
punctiform continuum (see \cite{23d}).  To Leibniz, points play only
the role of endpoints of line segments.  Thus the rejected
counterparts of the contradictory infinite wholes are not
infinitesimals but rather \emph{points} viewed as the simplest
constituents of a continuum.%
\footnote{The reciprocal relationship between Maxima and Minima can be
  formalized in modern mathematics as follows.  The unbounded infinity
  represented by the real line~$\R$ can be viewed as an increasing
  union of segments: ~$\R=\bigcup_{n\in\N}[-n,n]$.  In the reciprocal
  picture, we have the decreasing
  intersection~$\bigcap_{n\in\N}\big[-\frac1n, \frac1n\big]$ which is
  a single point (the origin).  Thus, the counterparts of infinite
  wholes are points (not infinitesimals).}

It is therefore problematic to assimilate -- as Rabouin and Arthur do
in \cite{Ra20} -- infinitesimals to infinite wholes in the matter of
inconsistency.  In 2022, Arthur again claims that
\begin{enumerate}\item[]
These unassignables, however, cannot be understood as actual
infinitesimals, since the notion of the actually infinitely small
contains a contradiction.  \cite[p.\;320]{Ar22}
\end{enumerate}
Arthur has repeatedly misinterpreted the contradiction involved.
According to Leibniz, the contradiction inherent in the notion of an
infinite whole, i.e., the \emph{infinitum interminatum}, does not
affect the \emph{infinita terminata} and their reciprocals, the
infinitesimals.  The most significant difference is that the
\emph{infinita terminata} are magnitudes, whereas infinite wholes are
multitudes.  While Leibniz considered the latter to be contrary to the
part-whole principle and therefore contradictory, infinite wholes will
be formalized by Cantor under the name of cardinalities, or
transfinite numbers.  By assimilating infinitesimals to infinite
wholes, Arthur and others are in essence attempting to invert a
cardinality so as to obtain an infinitesimal -- but to Leibniz, the
contradictory counterparts of infinite wholes are (not infinitesimals
but) \emph{points} viewed as constituent parts of the continuum.%
\footnote{Some historians are aware of the distinction and avoid
  committing Ishiguro and Arthur's error of conflating infinite number
  and infinite whole.  Thus, Spalt wrote: ``Johann Bernoulli
  {\ldots}~knew his friend Leibniz to be a philosophical thinker and,
  as such convinced that infinite `wholes' do not exist.  Therefore,
  he had to be careful in talking to Leibniz about `infinitely large'
  numbers'' \cite[p.\;53]{Sp22}.}

\subsection{Fraenkel shocked by inversion of cardinalities}
\label{s15}

Reacting to contemporary attempts to define infinitesimals, Abraham
Fraenkel wrote:
\begin{enumerate}\item[] 
I was deeply shocked at how infinity was treated in the Marburg
school, {\ldots}~wherein the infinitesimal is brought into direct
correspondence with Georg Cantor's transfinite numbers.
\cite[p.\;85]{Fr16}
\end{enumerate}
Here Fraenkel was referring to the work of Paul Natorp (1854--1924).
Fraenkel went on to describe Natorp's shocking treatment as follows:
\begin{enumerate}\item[]
If following Cantor an infinitely large `number' is denoted as~$w$,
and if the ratio~$1:w=x:1$ is formulated, then~$x$ must be an
`infinitely small number'.  (ad loc., note~18)
\end{enumerate}
Accordingly, Natorp confused magnitudes and multitudes, i.e.,
transfinite numbers.  What Fraenkel pointed out is that the latter
cannot be inverted to obtain an infinitesimal.  Fraenkel would have
likely been just as shocked to discover that, a century later, a
similar confusion persists in the writings of Arthur (see
Section~\ref{s14}) and other Leibniz scholars in the Ishiguro school;
see further in \cite{18i}.

Both in Leibniz and Cauchy scholarship, preoccupation with the
distinction between potential infinity and an infinite whole,
sometimes called actual infinity, is often a way of changing the
subject so as to deny that they ever used genuine infinitesimals.
Thus, Boyer writes: 
\begin{enumerate}\item[]
Cauchy and Weierstrass saw only paradox in attempts to identify an
actual or `completed' infinity in mathematics, believing that the
infinitely large and small indicated nothing more than the
potentiality of Aristotle.  \cite[p.\;612]{Bo68}
\end{enumerate}
What is involved is a conflation of two senses of the adjective
\emph{actual}:

\begin{enumerate}
\item[(A)] as in \emph{actual infinity} (or infinite whole) versus
  \emph{potential infinity}; and
\item[(B)] as in \emph{actual infinitesimal} number or quantity (i.e.,
  a genuine infinitesimal in the sense of violating Euclid's
  Definition V.4 when compared to~$1$), as opposed to a smaller and
  smaller ``ordinary'' number.  
\end{enumerate}
The term \emph{actual infinity} as used in set theory can refer to an
infinite multitude taken as a whole, whereas infinite number in the
sense of Leibniz can also refer to magnitude or quantity, not
multitude; for details see \cite{23d}.
%
%
The distinction between potential and actual infinity is distinct and
independent from the question whether Leibniz and Cauchy used genuine
infinitesimals.  We will examine the modern debate over Cauchy's
infinitesimals in Sections~\ref{s2} and \ref{s3}.

\section{Debate over Cauchy}
\label{s2}

\epigraph{Sim\'eon-Denis Poisson (1781--1840) {\ldots}~was able to
  promote [the \emph{infi\-niment petits}] to far-reaching
  dissemination and impact because of his central position within the
  French educational system.  ---Schubring~(2005)}

\subsection
{Variable quantities, infinitesimals and limits in the \emph{Cours}}
\label{s21b}

In his \emph{Cours d'Analyse} \cite{Ca21}, Cauchy laid foundations for
analysis that were characterized by the following features:
\begin{enumerate}
\item
Cauchy did not give an~$\epsilon$-$\delta$ definition of limit;
\item
Cauchy did not define the notion of continuity of a function in terms
of limits;
\item
Cauchy's notion of limit is similar to that of his teacher Lacroix
(who is not generally thought of as a pioneer of the Weierstrassian
\emph{Epsilontik});
\item
Cauchy's final definition of continuity of a function~$f$, emphasized
in italics in his \emph{Cours d'Analyse}, stipulates that an
infinitesimal increment~$\alpha$ must always produce an infinitesimal
change \mbox{$f(x+\alpha)-f(x)$} in the function.
\end{enumerate}
See further in \cite{20a} and \cite{21f}.  There have been varying
interpretations of what Cauchy meant by the term \emph{infinitesimal}.
He gave the following definition of infinitesimals:
\begin{enumerate}\item[]
When the successive numerical values of such a variable decrease
indefinitely, in such a way as to fall below any given number, this
variable becomes what we call \emph{infinitesimal}, or an
\emph{infinitely small quantity}.  A variable of this kind has zero as
its limit.%
\footnote{``Lorsque les valeurs num\'eriques successives d'une m\^eme
  variable d\'ecroissent ind\'efiniment, de mani\`ere \`a s'abaisser
  au-dessous de tout nombre donn\'e, cette variable devient ce qu'on
  nomme un \emph{infiniment petit} ou une quantit\'e \emph{infiniment
    petite}.  Une variable de cette esp\`ece a z\'ero pour limite''
  Cauchy \cite[p.\;4]{Ca21}.}
\cite[p.\;7]{BS}
\end{enumerate}
Interpretations have varied with regard to Cauchy's term
\emph{devient} (becomes).  Grabiner \cite{Gr81} holds that Cauchy's
\emph{becomes} is equivalent to \emph{is}, so that an infinitesimal is
nothing but a variable quantity tending to zero, ruling out any
associated non-Archimedean phenomena.  Others have argued that the
term \emph{becomes} implies a process involving a change of nature%
\footnote{
Such a process can be formalized in modern mathematics in terms of a
suitable equivalence relation; see e.g., \cite{21d}.}
and have, accordingly, interpreted Cauchy's infinitesimals in
non-Archimedean terms.  Regardless of the meaning of Cauchy's
definition of infinitesimals in the \emph{Cours}, he used genuine
infinitesimals in his later books and research articles, as documented
in \cite{20b} (see below).

Detlef Laugwitz, following Robinson \cite[pp.\;269--276]{Ro66},
proposed an interpretation of Cauchyan analysis in a series of
articles starting in the 1980s (see e.g., \cite{La87} and
\cite{La89}), sparking a historical debate that is still current.

In 2020, Bair et al.\;\cite{20b} explored several applications that
Cauchy made of infinitesimals in his work beyond his \emph{Ecole
  Polytechnique} textbooks, in fields ranging from centers of
curvature (see Section~\ref{s24b} below) to convergence of series of
functions (see Section~\ref{s36}) to integral geometry (see
Section~\ref{s29}).

In 2022, Gert Schubring (GS) published a lengthy review \cite{Sc22} of
Bair et al.~\cite{20b} for MathSciNet.  Such attention is surely
appreciated by every researcher (even though the review is unrefereed,
i.e., non-peer-reviewed); hopefully it can mark the beginning of a
meaningful dialog or informed debate.

It turns out, however, that many of GS's comments in the review are at
odds with what he wrote in his 2005 book \cite{Sc05}, which he would
go on to describe in \cite[p.\;527]{Sc16} as his ``key publication.''
We will document multiple shifts in GS's position between 2005 and
2022, including his book's acknowledgment of Poisson's broad influence
in promoting Poisson's version of infinitesimals, which rivaled the
type of \emph{infinitesimal lore} that GS chose exclusively to
emphasize in his review.  We will also compare GS's reactions to the
Cauchy scholarship of Laugwitz and that of Bair et al.  The
philosophical assumptions underpinning GS's position are analyzed in
Section~\ref{s39}.

\subsection{Members of a group}
\label{s21}

In his opening remarks, GS makes it clear that he is targeting
not merely the article under review, but an entire program of
re-evaluation of the history of analysis pursued by a large group of
scholars:
\begin{enumerate}\item[]
This paper is written by members of a group that for several years has
been leading a \emph{crusade against the historiography} of
mathematics.~{\ldots}~the group seems to consist of at least 22
mathematicians and philosophers, {\ldots} \cite[p.\,1]{Sc22}%
\footnote{Page numbers here and below refer to the pdf of GS's
  review.}
(emphasis on ``crusade against the historiography'' added)
\end{enumerate}
With regard to GS's choice of wording (``crusade against, etc.'') in
describing the work of scholars he happens to disagree with, it is
worth recalling his track record of colorful terminology targeting the
work of Laugwitz, as for example the following sarcastic comment:
\begin{enumerate}\item[]
``[Giusti's 1984 article] spurred Laugwitz to even more detailed
  attempts to banish the error and confirm that Cauchy had used
  hyper-real numbers.  {\ldots} (see Laugwitz 1990, 21).''  (Schubring
  \cite[p.\;432]{Sc05})%
\footnote{\label{f4}For the record, we note that Laugwitz never
  attributed hyperreal numbers to Cauchy, either in his 1990 article
  \cite{La90} cited by GS or anywhere else.  For an analysis of GS's
  misrepresentation of Laugwitz, see \cite[Section~6.1]{17e} and
  \cite[Section~4.5, pp.\;278--279]{18e}.  See further in
  Section~\ref{s97} on Laugwitz's take on procedures.}
\end{enumerate}
Further epithets are sampled in Section~\ref{s212}.  One of the works
GS quotes in 2022 is the article Bair et al.~\cite{13a} going back to
2013, indicating that he is targeting at least a decade of work by the
22 scholars he mentioned.  GS explains that
\begin{enumerate}\item[]
To enable the reader to understand the issues at stake and to situate
their development, this review is a bit more extended than usual.
\cite[p.\,1]{Sc22}
\end{enumerate}
We will examine GS's take on ``the issues at stake,'' evaluate how
successful his attempts to ``situate their development'' are, and
compare his views as expressed in 2005 and in 2022.

\subsection{Far-reaching dissemination: Carnot, Coriolis, Poisson}
\label{f5b}

Concerning the historical period of Cauchy's activity, GS claims the
following in 2022:
\begin{enumerate}\item[]
Nobody so far, including Cauchy himself, had thought of the calculus
in terms of a non-Archimedean continuum.  \cite[p.\;2]{Sc22}
\end{enumerate}
%
%
However, an examination of GS's key publication \cite{Sc05} indicates
that he did not always feel this way.  In 2005, he acknowledged the
following:
\begin{enumerate}\item[]
[T]here was a propagation of the \emph{infiniment petits} quite
deviant [sic] from Cauchy’s conceptualization in France, an attempt to
popularize them as \emph{actually} infinitely small quantities.  This
version would hence have had a claim on \emph{anticipating
  non-standard analysis} if its propagators had been able to know that
the latter would some day exist.  \cite[p.\;575]{Sc05} (emphasis on
``actually'' in the original; emphasis on ``anticipating non-standard
analysis'' added)
\end{enumerate}
GS goes on to identify the propagator:
\begin{enumerate}\item[]
The conceptions were those of Sim\'eon-Denis Poisson (1781--1840), who
was able to promote them to \emph{far-reaching dissemination} and
impact because of his central position within the French educational
system.  (ibid.; emphasis on ``far-reaching dissemination'' added)
\end{enumerate}
As member of the \emph{Conseil royal de l'instruction publique},
Poisson succeeded in promoting infinitesimals to officially prescribed
status at the national level in 1837:
\begin{enumerate}\item[]
Poisson prevailed in having his favorite method of using the
\emph{infiniment petits} centrally prescribed as compulsory of all the
\emph{coll\`eges}: ``The two geometry lessons will remain appended to
the \emph{troisi\`eme} class; but this teaching will be based on the
method of infinitely small quantities.''  (Schubring
\cite[p.\;585]{Sc05} based on Belhoste \cite[p.\,147]{Be95})
\end{enumerate}
A decree issued during the following year was even more specific:
\begin{enumerate}\item[]
A subsequent decree dated October 9, 1838, specified, by presenting a
first detailed mathematics curriculum, how geometry was to be taught
on the basis of the \emph{infiniment petits}:
\begin{enumerate}\item[]
In plane geometry for the \emph{troisi\`eme}, curves were to be
conceived of as polygons having an infinite number of sides; in
particular, circles were to be conceived of as regular polygons.
\end{enumerate}
(Schubring \cite[p.\;585]{Sc05} based on Belhoste
\cite[p.\,148]{Be95})
\end{enumerate}
Such a conception of a curve as an infinite-sided polygon, prescribed
for the French mathematics curriculum in 1838, goes back to Leibniz
and even earlier (a point not mentioned by GS); see further in
\cite[Section~2.3]{21a}.

Nor was Poisson the only first-rate mathematician to promote the
dissemination of genuine infinitesimals.  As noted by
Grattan-Guinness, Coriolis wrote the following:
\begin{enumerate}\item[]
``The approbation of the \emph{Conseil Royal} of the
  \emph{Universit\'e} will be equal to giving an appropriate direction
  to these works and to establishing everywhere the same language
  founded upon the infinitely small.''  (Coriolis as translated by
  Grattan-Guinness in \cite[p.\,1262]{Gr90})
\end{enumerate}
Similar remarks apply to Carnot.%
\footnote{See note~\ref{f22b}.}

In short, GS's key publication acknowledged a ``far-reaching
dissemination and impact'' of genuine infinitesimals in France, which
was moreover endorsed at the national level by the \emph{Conseil
  Royal}.%
\footnote{\label{f20}In 2022, GS claims that ``[Bair et al.] call
  [Poisson], wrongly, a member of the CP [\emph{Conseil de
      Perfectionnement} of the Ecole Polytechnique], while he was in
  fact an omnipotent member of the ministry's \emph{Conseil royal
    d'instruction publique}'' \cite[p.\;6]{Sc22}.  GS's claim is in
  error.  Indeed, Gilain mentions the ``examinateurs de
  math\'e\-matiques, \emph{Poisson} et \emph{de Prony}, qui animaient
  en g\'en\'eral la commission programme du CP'' \cite[\S 32]{Gi89}.
  GS's error is particularly surprising given his emphasis on
  Poisson's perceived vice of the \emph{cumul des mandats} in
  \cite[p.\;576]{Sc05}, so that his being on the \emph{Conseil Royal}
  is not inconsistent with the fact that he ``championed the use of
  infinitesimals through [his] influence on the \emph{Conseil de
    Perfectionnement} (CP)'' as mentioned by Bair et
  al.~\cite[p.\,142]{20b}.}
Very little of the above information concerning the extent of such
dissemination trickled down to his 2022 review.  To deny such a 19th
century lore of genuine infinitesimals is ``to wrench Cauchy's ideas
out of their historical context'' (cf. \cite{Ar22}).

Specifically, GS's 2005 talk of 19th century anticipation of
nonstandard analysis (NSA) undermines his 2022 claim that ``nobody
{\ldots}\;had thought of the calculus in terms of a non-Archimedean
continuum.''  If GS changed his mind about the ``far-reaching
dissemination'' of Poisson's non-Archimedean approach, he did not tell
his readers about such a change of heart.%
\footnote{\label{f22b}In 2022, GS writes: ``Lazare Carnot devoted
  himself, in his pre-infinitesimal periods, to elaborating the
  concept of null sequences for variables with limit zero''
  \cite[p.\;2]{Sc22}.  If there was a pre-infinitesimal period, there
  must have been a subsequent infinitesimal period, as well.  GS seems
  to acknowledge implicitly that genuine infinitesimals (or the
  \emph{pseudo-infiniment petits} in GS's parlance; see
  Section~\ref{s24}) were practiced by Carnot, as well.  For an
  analysis of Carnot's conception in relation to Leibniz's, see
  Barreau \cite{Ba89}.}
The 2022 claim is hardly compatible with GS's stated
goal of properly situating the development of infinitesimals in their
historical context (see Section~\ref{s21}).

\subsection{Epsilons, small and infinitesimal}
\label{s24b}

In a 1826 work on differential geometry, Cauchy develops a formula for
the radius of curvature~$\rho$ of a plane curve parametrized by
arclength~$s$.  Such a formula is equivalent to the
formula~$\frac1\rho=\frac{d\tau}{ds}$ in modern notation, where~$\tau$
is similar to the modern polar coordinate angle~$\theta$ of the
tangent vector to the curve.

Cauchy's proof of the formula was analyzed in \cite[Section~5]{20b}
(for a summary see below, Section~\ref{s211}, item~5).  The conclusion
was that the symbol~$\varepsilon$ \emph{as used there} denotes a
genuine infinitesimal on par with the quantities~$\Delta\tau$,~$\Delta
s$,~$\Delta x$, and~$\Delta y$ also used there.
%
%

GS rejects such a conclusion in the following terms:
\begin{enumerate}\item[]
[T]he use of~$\varepsilon$ is known as being due to Cauchy himself,
and this occurred apparently for the first time in his 1823 textbook
on the differential calculus, in which he used it for the derivative,
demanding `Let~$\delta,\varepsilon$ be two very small numbers',
clearly understanding them to have finite values.  \cite[p.\;2]{Sc22}
\end{enumerate}
The problem with GS's claim is two-fold.
\begin{enumerate}
\item
GS's remark with regard to the 1823 textbook is in error, as Cauchy
had already used ~$\varepsilon$ in this sense in his 1821 textbook
\cite[Section~2.3, Theorem~1]{Ca21}; see further in
\cite[Section~3.8]{19a} (see also Section~\ref{f18} below for an
analysis of GS's oversight).
\item
Cauchy used the notation~$\varepsilon$ on more than one occasion.  But
if ``infinitesimal'' and ``very small'' meant the same thing to
Cauchy, why did he use different terms, especially given the existence
of a widely disseminated Poissonian infinitesimal tradition (see
Section~\ref{f5b})?  
\end{enumerate}
Here GS seems to have overlooked the possibility that~$\varepsilon$
may \emph{not} have had the same meaning in Cauchy's 1826 text on
differential geometry as in his 1821 and 1823 calculus textbooks.

\subsection{Schubring \emph{vs} Grabiner}
\label{f18}

GS's oversight mentioned in Section~\ref{s24b} is all the more
puzzling since this particular 1821 occurrence of~$\varepsilon$ is one
of the key pieces of evidence used by Grabiner to argue Cauchy's
pioneering work on ``manipulating algebraic inequalities.''  Grabiner
translates the \emph{Cours d'Analyse} as follows:
\begin{enumerate}\item[]
``Designate by~$\varepsilon$ a number as small as desired.  Since the
  increasing values of~$x$ will make the difference~$f(x+1)-f(x)$
  converge to the limit~$k$, we can give to~$h$ a value sufficiently
  large so that,~$x$ being equal to or greater than~$h$, the
  difference in question is included between~$k-\varepsilon$
  and~$k+\varepsilon$.''  (Cauchy as translated by Grabiner in
  \cite[p.\;8]{Gr81})
\end{enumerate}
Grabiner goes on to conclude: ``This is hard to improve on'' (ibid.),
indicating that she views Cauchy's calculation as a convincing example
of an~$\varepsilon,\delta$ argument (a point already mentioned by
Freudenthal \cite[p.\;137]{Fr81}).  In 2016, GS sharply disagreed with
Grabiner in the following terms: 
\begin{enumerate}\item[]
I am criticizing historiographical approaches like that of Judith
Grabiner where one sees epsilon-delta already realized in Cauchy.
\cite[p.\;530]{Sc16}
\end{enumerate}
It emerges that Schubring criticized Grabiner without properly
examining the evidence she had presented.

\subsection{Non-standard numbers}
\label{s23b}

While in his key publication GS acknowledged the existence of a widely
disseminated non-Archimedean conception of infinitesimal and infinite
numbers as advocated by Poisson (see Section~\ref{f5b}), by 2022 we
find GS claiming the following:
\begin{enumerate}\item[]
Cauchy used for them the term of traditional \emph{lore} to speak of
arbitrarily large numbers --- as mathematicians did throughout the
18th century, and likewise used by Weierstra{\ss} [see, for instance,
  K.\;Viertel {\ldots}\;--- none of them ever thinking of
  \emph{non-standard numbers}.  \cite[p.\;3]{Sc22} (emphasis on
  ``lore'' and ``non-standard numbers'' added)
\end{enumerate}
GS's reference to Weierstrass is particularly revealing.  If there did
exist a widely disseminated non-Archimedean tradition of genuine
infinitesimals in Cauchy's time (as GS acknowledged in his key
publication), why should we automatically assume that when
mathematicians used the terms \emph{infinitely small} or
\emph{infinitely large}, they necessarily meant it in the
Weierstrassian sense, as shorthand for more long-winded
non-infinitesimal arguments?

While GS emphasizes his student Viertel's work,
%
%
%
Weierstrass' use of the term \emph{infinitesimal} in
\cite[p.\;74]{We86} in a figurative sense was quoted in Bascelli et
al.\;\cite[Section~2.1]{16a}, a fact not acknowledged by~GS.  The
existence of such traditional lore was mentioned by Laugwitz in 1989
in the following terms:
\begin{enumerate}\item[]
In 1815 [Cauchy] does not use infinitesimals but only `very small
numbers', in a naive and pragmatic manner.  \cite[p.\;232]{La89}
\end{enumerate}
Laugwitz goes on to argue that the change occurred in the early 1820s,
when Cauchy started using genuine infinitesimals.

Appealing to a ``traditional lore'' as evidence that Cauchy's infinite
numbers were merely large ``ordinary'' numbers is a non-sequitur,
given the dual nature of the said lore, as we elaborate further in
Section~\ref{s24}.

\subsection{Dual nature of 19th century  lore}  
\label{s24}

Even though the term \emph{lore} was absent from the analysis of the
18th and 19th centuries in his key publication from 2005, by 2022 GS
speaks of infinitesimal and infinite numbers as parts of a
``traditional lore to speak of arbitrarily large numbers --- as
mathematicians did throughout the 18th century'' (see
Section~\ref{s23b}).  Accordingly, he describes Cauchy's phrase
\begin{quote}
``si l'on d\'esigne par~$\varepsilon$ un nombre infiniment petit'' 
\end{quote}
(if one denotes by $\varepsilon$ an infinitely small number) as
``standard lore for expressing an arbitrarily small number''
\cite[p.\;3]{Sc22}.

What was the nature of the said \emph{lore}?  One significant
development that GS failed to mention is the 19th century debate about
infinitesimals in Belgium and Luxembourg.  Jean-Nicolas No\"el
(1783--1867) at the University of Li\`ege and Jean Joseph Manilius
(1807--1869) at the University of Ghent (Gand) \cite[p.\;8]{Bo66} were
advocates of the use of genuine infinitesimals.  They introduced both
the Leibnizian distinction between assignable and inassign\-able
quantities, and Leibniz's definition of infinitesimal as smaller than
any assignable quantity (see Section~\ref{s13}).  Their opponents were
led by Ernest Lamarle (1806--1875) similarly at the University of
Gand; see further in \cite{Bo66}, \cite{Ba19}, and \cite{23c}.  It
emerges that there were distinct and \emph{rival} infinitesimal lores
in both France and Belgium at the time.
%
%
%

In 2022, GS does mention a criticism of genuine infinitesimals ``in
Belgium in 1887 by Paul Mansion'' \cite[p.\;7]{Sc22}, but fails to
clarify its context.  The context was indeed the debate opposing
No\"el/Manilius and Lamarle; Mansion \cite{Ma87} (at the University of
Gand like Manilius and Lamarle) was presenting a rebuttal of the
position held by No\"el and Manilius, a fact GS fails to mention.

The Belgian debate is significant in the context of GS's references to
an infinitesimal \emph{lore}, which he claims to amount to viewing an
infinitesimal as an arbitrarily small number.  There were surely
infinitesimal lores in France, Belgium, and elsewhere as he claims,
but GS appears to be only reporting one of them in 2022.  In sum, by
seeking to portray a monolithic lore, GS fails to clarify the dual
nature of the infinitesimal lore.%
\footnote{It is therefore ironic that GS should accuse the authors of
  \cite{20b} of holding such monolithic views, when he writes: ``The
  authors are so enraptured by their conviction that Cauchy conceived
  of numbers within a non-Archimedean continuum that they consider
  historical mathematics anachronistically only in terms of this
  continuum'' \cite[p.\;2]{Sc22}.  Actually, it is GS who painted a
  reductionist monolithic picture of infinitesimal lore in his review,
  whereas the authors of \cite{20b} recognize the dual nature
  thereof.}

GS was more forthcoming in 2005 with details on such a rival approach
to infinitesimals, and even gave it a name: \emph{pseudo-infiniment
  petits} (pseudo-infinitely small quantities), adopting Mansion's
terminology:
\begin{enumerate}\item[]
Mansion called quantities thus conceived of as \emph{pseudo-infiniment
  petits}, with the intention of recalling the \emph{obvious
  contradiction} this \emph{definition} contained.
%
%
\cite[p.\;583]{Sc05} (emphasis on ``obvious contradiction'' and
``definition'' added)
\end{enumerate}
The definition in question was formulated as follows by Mansion:
\begin{enumerate}\item[]
Certain geometers have given yet another meaning to the word
\emph{infinitely small}.  According to them, there exist quantities
different from zero which are yet smaller than every assignable
magnitude.%
\footnote{``Certains g\'eom\`etres ont donn\'e un autre sens encore au
  mot \emph{infiniment petit}.  D'apr\`es eux, il existe \emph{des
    quantit\'es diff\'erentes de z\'ero et qui sont cependant
    inf\'erieures \`a toute grandeur assignable}''
  \cite[p.\;214]{Ma87}.}
\end{enumerate}
The definition criticized by Mansion is essentially the Leibnizian
definition of infinitesimals as smaller than every assignable quantity
(see Section~\ref{s13}).%
\footnote{\label{f7}The fact that this definition of infinitesimals
  admits a modern formalisation (see note~\ref{f3}) provides evidence
  against the Mansion--Schubring opinion that the definition involved
  an ``obvious contradiction.''  Mansion's opinion is perhaps
  understandable as it may have been common in the 1880s, but
  Schubring's is less~so.}
GS's key publication devotes a number of pages to the subject in
\cite[pp.\;583--593]{Sc05} where the term \emph{pseudo-infiniment
  petits} is mentioned six times.  Alas, very little of such
pseudo-infinitesimal lore trickled down to his 2022 text, which
therefore fails to situate the development of infinitesimals in their
proper historical context with all its complexity.

\subsection{What is a number, ontologically?}
\label{f3b}

The \emph{Cours d'Analyse} includes a lengthy appendix called
\emph{Note I} \cite[pp.\;403--437]{Ca21}.  GS claims that Cauchy only
recognized finite numbers in his \emph{Note I}:
\begin{enumerate}\item[]
{\ldots} the authors still refuse to take notice of Cauchy's very
explicit discussion of what the number system means for him, in the
introduction and in the extensive \emph{Note I} in his \emph{Cours
  d'analyse} of 1821 {\ldots}~The status of ``number'' is attributed
exclusively to positive integer numbers, called by Cauchy absolute
numbers.  {\ldots}~In their crusade to capture Cauchy as a
non-standard analysis forerunner, they systematically ignore Cauchy's
own explicit affirmations that for him only \emph{finite numbers} are
admitted as numbers.  \cite[p.\;3]{Sc22} (emphasis on ``Note I'' in
the original; emphasis on ``finite numbers'' added)
\end{enumerate}
The problem with GS's claim is four-fold.

\begin{enumerate}
\item
A reading of Cauchy's \emph{Note I} reveals not only that he never
claimed that all numbers are finite, but on the contrary that he
envisioned the possibility that they may not be.  Indeed, following
the introduction of exponentiation~$A^B$, Cauchy presents the
following formula and comment:
\[
A^0=1.
\]
\begin{quote}
We assume however that the value of the number $A$ remains finite and
differs from zero.%
\footnote{``Nous supposons toutefois que la valeur du nombre~$A$ reste
  finie et diff\`ere de ze\'ro'' \cite[p.\;416]{Ca21}.}
\end{quote}
Thus, Cauchy finds it necessary to stipulate the \emph{condition} that
the number~$A$ should be finite (possibly because he is thinking of
indeterminate forms of type~$\infty^0$) to ensure the validity of the
formula.
\item
GS's claim that in \emph{Note I} all numbers are positive integers is
contradicted by Cauchy's analysis of both rational and irrational
numbers in \emph{Note I} \cite[p.\;409]{Ca21}.
\item
GS claims that
\begin{enumerate}\item[]
Cauchy followed Carnot%
\footnote{Schubring's position on Carnot is obscure; see
  note~\ref{f22b}.}
in interpreting numbers ontologically and admitting only absolute
numbers.  \cite[p.\;3]{Sc22}
\end{enumerate}
However, Cauchy's ``absolute numbers'' are simply \emph{unsigned}
numbers, as is evident from Cauchy's Turin lectures where they are
referred to as \emph{numeri assoluti} \cite[p.\,152]{Ca30}.  Little
can be derived from Cauchy's comments on unsigned numbers as regards
the \emph{ontology} of his numbers.
\item
In 2022, GS discusses the passage from the text on the centers of
curvature where Cauchy writes ``si l'on d\'esigne par~$\varepsilon$ un
nombre infiniment petit.''  GS's claim that Cauchyan numbers are
necessarily positive integers is at odds with Cauchy's reference to an
infinitesimal as a \emph{number}, a reference GS himself quotes.
\end{enumerate}

What Cauchy did affirm, similarly to Leibniz (see Section~\ref{s13}),
was that infinite wholes are contradictory, being contrary to the
part-whole principle.  Thus, in the \emph{Sept Le\c{c}ons}, Cauchy
summarizes Galileo's paradox (comparing numbers and their squares),
and concludes: 
\begin{enumerate}\item[]
The proof that we just recalled was given for the first time by
Galileo.%
\footnote{``[L]a d\'emonstration que nous venons de rappeler a \'et\'e
  donn\'ee pour la premi\`ere fois par Galil\'ee''
  \cite[p.\;422]{Ca33}.}
\end{enumerate}
%
%
GS's claim that ``[Cauchy] says that one can prove mathematically that
the assumption of a number `infinite' would lead to manifest
contradictions'' \cite[p.\;449]{Sc05} equivocates on the meaning of
the term \emph{number}: Cauchy was referring only to the impossibility
of \emph{infinite wholes}.%
\footnote{\label{f19}Leibniz's position was similar; see main text at
  note~\ref{f8}.}
See further in Laugwitz \cite[p.\;201]{La89}.

In his key publication, GS acknowledged \cite[pp.\;445,
  448]{Sc05} that Abbot Moigno edited and published Cauchy's
\emph{Sept Le\c{c}ons de Physique G\'en\'erale} already \emph{after}
Cauchy's death.  Yet in 2022 GS is willing to rely upon the
good Abbot
%
%
to represent Cauchy's views faithfully when Moigno claims in an
Appendix (written by himself) to \emph{Sept Le\c{c}ons} that ``A
number being actually infinite is impossible; every number is
essentially finite'' \cite[p.\;4]{Sc22}.%
\footnote{%
%
As noted by Redondi, ``Parmi tous les myst\`eres de la raison, Moigno
retient celui de l'impossibilit\'e logique d'un nombre actuellement
infini'' \cite[p.\;217]{Re88}.}
Moigno believed this, but did Cauchy (see Section~\ref{s36})?  GS has
much to say about Moigno's and Mansion's criticism of Poisson's
genuine infinitesimals,%
\footnote{I.e., the \emph{pseudo-infiniment petits} in GS's parlance;
  see Section~\ref{s24}.}
but there is one name conspicuously absent from GS's list of the
critics of Poisson: Cauchy himself.  Significantly, Moigno's lengthy
introduction was deleted when the \emph{Sept Le\c{c}ons} were included
in the \emph{Oeuvres Compl\`etes} under the direction of the French
Academy of Sciences; see \cite[p.\;412]{Ca33}.  Similarly left out was
Moigno's lengthy appendix pompously entitled
\begin{quote}
``Sur l'impossibilit\'e du nombre actuellement infini, l'antiquit\'e
  de l'homme, la science dans ses rapports avec la foi'' (ibid.).
\end{quote}

GS's assumption that Moigno and Mansion expressed Cauchy's view
remains in the realm of opinion rather than a supported position.
Describing Moigno as Cauchy's ``alter religious-philosophical ego'' as
GS does in \cite[p.\;3]{Sc22} carries little persuasive force.  It
must be noted that already in his key publication, GS expressed his
appreciation of the good Abbot: 
\begin{enumerate}\item[]
Moigno not only rejects Poisson on the basis of the concepts shared at
that time; he simultaneously discloses the \emph{massive
  contradiction} in the French mathematical community.
\cite[p.\;455]{Sc05} (emphasis added)
\end{enumerate}
GS's bold attribution of a ``massive contradiction'' to the ``French
mathematical community'' is not accepted by all historians.%
\footnote{For an analysis of the Moigno--Schubring ``massive
  contradiction'' see note~\ref{f7}.}

\section{Sum theorem, integral geometry, and continuity}
\label{s3}

\subsection{Sum theorem}
\label{s36}

Cauchy's sum theorem concerning conditions for continuity of the sum
of a series of continuous functions has long been the subject of a
controversy.  In 1821, Cauchy published a version of the theorem in
his \emph{Cours d'Analyse} \cite{Ca21}.  Abel and others eventually
pointed out that the theorem seems to ``suffer exceptions.''  In 1853,
Cauchy published a version of the theorem with an apparently modified
hypothesis in the article \cite{Ca53}.  The issue is summarized in
Section~\ref{s33b}.  Some scholars have argued that the modified
hypothesis is equivalent to uniform convergence.  Robinson's
interpretation includes the following two items:
\begin{enumerate}
\item
\label{l1}
the hypothesis of the 1821 theorem required only convergence at
standard points, and therefore without additional assumptions, the
theorem was incorrect as stated \cite[pp.\;271--272]{Ro66};
\item
\label{l2}
the hypothesis of the 1853 theorem required convergence at all points,
including infinitesimals, resulting in a correct theorem when
interpreted in nonstandard analysis \cite[p.\;273]{Ro66}.%
\footnote{\label{f22}See Section~\ref{spalt} on Spalt's coverage of
  Robinson's interpretation.}
\end{enumerate}

The sum theorem was mentioned by Bair et al.~in
\cite[Section\;2]{20b}.  Reacting to this mention, Schubring claims
the following:

\begin{enumerate}\item[] 
[T]he authors' principal aim here is to show that in these
applications, Cauchy used infinitesimals as numbers, hence in the
non-standard meaning.  Yet, right at the beginning, this intention
leads the authors to \emph{falsify} a text by Cauchy.  In their
initial section on summation of series (Bair et al.~2020, p.\,130),
they give a truncated quotation: ``He then states his convergence
theorem modulo a hypothesis that the sum~$u_n+u_{n+1}+\ldots+u_{n'-1}$
should be \emph{toujours infiniment petite pour des valeurs infiniment
  grandes des nombres entiers~$n$ et~$n' > n$ \ldots} (Cauchy [18],
1853, p.\;457)'' {\ldots} Cauchy's text had been, however:
\emph{{\ldots} la somme devient toujours infiniment petite, quand
  {\ldots}\,.}  Cauchy had, hence, unequivocally expressed that the
sum becomes as small as one wishes, thus dealing with a limit process.
\cite[pp.\;2--3]{Sc22} (emphasis of ``falsify'' added)
\end{enumerate}
In what way has Cauchy's text been allegedly falsified?  Consider the
following two items.

\begin{enumerate}
\item
\label{one}
Note the difference in wording: \cite{20b} quoted the Cauchyan passage
as saying ``toujours infiniment petite pour etc.''  whereas GS
quotes it as saying ``devient toujours infiniment petite, quand etc.''
Thus allegedly Cauchy's \emph{quand} was replaced by \emph{pour}.
\item
\label{two}
GS stresses that the Cauchyan passage was truncated through the
deletion of the verb \emph{devient}, which, as he claims, represents a
limiting process.
\end{enumerate}

As far as item \eqref{one} is concerned, note that on page 457 Cauchy
writes 
\begin{enumerate}\item[]
\emph{pour des valeurs infiniment grandes des nombres entiers~$n$ et
  $n'>n$},
\end{enumerate}
using \emph{pour}, exactly as quoted in \cite{20b}.  On the next page
458 (\emph{not} cited in~\cite{20b}), there is a complex variable
version of the theorem where Cauchy writes
\begin{enumerate}\item[]
\emph{quand on attribue des valeurs infiniment grandes aux nombres
  entiers~$n$ et~$n'>n$},
\end{enumerate}
using \emph{quand} as quoted by GS.  It emerges that GS was merely
looking at the wrong page in Cauchy's paper when reviewing Bair et
al.~\cite{20b}.

As far as item \eqref{two} is concerned, the interpretation of
Cauchy's \emph{devient} was discussed in detail in B{\l}aszczyk et
al.~\cite{17e} and elsewhere.  GS's suggestion that the issue is being
ignored by the group of scholars in question (see Section~\ref{s21})
is therefore baseless and misleading.  Cauchy's \emph{devient} was
also discussed in a 2011 article in \emph{Perspectives on Science}
\cite{11b}, and specifically in reference to Cauchy's sum theorem; and
more recently in Bascelli et al.~(\cite{18e}, 2018).  The rival
interpretations are summarized in~\cite{21f}.  The substantive issue
with regard to item \eqref{two} was outlined in Section~\ref{s21b}.

If GS wishes to go with an Archimedean interpretation, he would have
to interpret \emph{becomes} as simply \emph{is}.  Then GS would find
Cauchy asserting that 
\begin{quote}
``the sum~$u_n+u_{n+1}+\ldots+u_{n'-1}$ \emph{is} always infinitely
  small for infinitely large values of the integers~$n$ et~$n' > n$.''
\end{quote}
This is a reasonable interpretation of the passage that ultimately
hinges on the meaning of the term \emph{infinitely small}; see
Section~\ref{s21b}.  However, GS's assumption that the word
\emph{devient} ``represents a limiting process'' in a necessarily
Archimedean context, remains in the realm of opinion rather than a
position supported by evidence.

\subsection{Procedures and foundations in Spalt}
\label{spalt}

In a 2022 book, Detlef Spalt claims that
\begin{enumerate}\item[]
[Robinson's] result was: Cauchy's theorem is correct if we add one of
the two additional assumptions: (a) the series is uniformly convergent
or (b) the family~$(s_n(x))_n$ of partial sums is equicontinuous in
the interval.  \cite[p.\;239]{Sp22}
\end{enumerate}
However, Spalt's claim is true only with regard to Robinson's
interpretation of the 1821 formulation of the sum theorem.  As
mentioned in Section~\ref{s36}, item~\eqref{l1}, Robinson assumes that
convergence was required only at standard $x$ in 1821.  Meanwhile, as
mentioned in item \eqref{l2}, Robinson goes on to address Cauchy's
1853 formulation, and states clearly that Cauchy's term ``toujours''
refers to the extension of the convergence condition from standard~$x$
to all~$x$ (including infinitesimal values):
\begin{enumerate}\item[]
If we interpret this theorem in the sense of Non-standard Analysis, so
that `infiniment petite’ is taken to mean `infinitesimal’ and
translate `toujours’ by `for all~$x$’ (and not only `for all
standard~$x$’), then the condition introduced by Cauchy
{\ldots}~amounts precisely to uniform convergence.
\cite[p.\;273]{Ro66}
\end{enumerate}
Spalt's misrepresentation of Robinson's Cauchy scholarship (see
further in Section~\ref{sprob}) is symptomatic of a deeper problem.  A
colorful instance of conflation of procedures and foundations occurs
in Spalt, who wrote:
\begin{enumerate}\item[]
Nonstandard-analysis needs to perform some acts of conceptual
acrobatics in order to construct ``hyper-real'' numbers in a
mathematically acceptable way.  Here Robinson's construction stands
out, but without studying (at least) one semester of modern logic, one
is not able to follow his construction.  It is, however, unlikely that
Cauchy should have anticipated such a stilted concept in 1821 in any
possible sense.  \cite[p.\,134]{Sp22}
\end{enumerate}
This ``acrobatics'' passage overlooks the fact that the issue is not
how hard Robinson's construction is, but rather how hard it is for
Spalt (and Schubring) to appreciate the distinction between procedures
and foundations; see further in Section~\ref{s97}.

\subsection{Spalt on Robinson}
\label{sprob}

Oddly, two decades earlier (in 2002), Spalt did recognize that in
Robinson's interpretation, the 1853 result was correct:
\begin{enumerate}\item[]
After citing Cauchy's later formulation of the sum theorem from 1853,
Robinson concludes that the formulation given there `amounts exactly
to the uniform convergence in agreement with (i) above'.%
\footnote{``Nachdem Robinson die sp\"atere Formulierung des
  Summensatzes durch Cauchy aus dem Jahr 1853 angef\"uhrt hat, kommt
  er zu dem Schluss, die dort gegebene Formulierung laufe `genau auf
  die gleichm\"a{\ss}ige Konvergenz in \"Ubereinstimmung mit (i) oben
  hinaus'{}'' \cite[p.\;297]{Sp02}.}
\end{enumerate}
In general, inaccuracies abound in Spalt's reporting on nonstandard
analysis.  Thus, he claimed that the intermediate value theorem is
false in nonstandard analysis: 
\begin{enumerate}\item[]
Take for example the Intermediate value theorem.  Every continuous
function which has two different values takes on each intermediate
value in between.  This is a basic theorem of classical analysis.  But
it is not true in non-standard analysis, {\ldots} \cite[p.\,167]{Sp01}
\end{enumerate}
Of course the intermediate value theorem is in the scope of the
transfer principle and therefore is as true in nonstandard analysis as
in classical analysis, and moreover can be proved using
infinitesimals; see e.g., \cite[p.\;67]{Ro66}.

\subsection{Integral geometry}
\label{s29}

One of Cauchy's goals in \cite{Ca26} is a formula in integral geometry
expressing the length of a curve in terms of an average of the lengths
of its projections to the pencil of lines through the origin.  There
are two aspects of the problem, described in \cite{20b} as ``the
curve" and ``the Grassmannian,'' i.e., circle of directions from a
point.  Cauchy subdivides the curve into infinitesimal subsegments
(\emph{\'el\'ements infiniment petits}).  By contrast, he approximates
the Grassmannian circle by a regular~$n$-gon, and studies the
asymptotic behavior of the resulting approximation as~$n$ tends to
infinity.  If an infinitesimal merely meant a sequence to Cauchy, then
there shouldn't be any difference in Cauchy's treatment of the curve
and the Grassmannian; both should be sequences.  The fact that Cauchy
does treat them differently suggests that his infinitesimals are not
merely sequences.

Concerning the analysis in \cite{20b} of Cauchy's theorem in integral
geometry, GS writes the following:

\begin{enumerate}\item[]
\emph{Lengths}. In a section on integral geometry, the authors refer
to a publication by Cauchy where he ``exploits a decomposition of a
curve into infinitesimal length elements (respectively, of a surface
into infinitesimal area elements)'' (Bair et al.~2020, p.\,130). The
key formulation for them is Cauchy's statement (Bair et al.~2020,
p.\,132):
\begin{enumerate}\item[]
Le th\'eor\`eme II \'etant ainsi d\'emontr\'e pour le cas particulier
o\`u la quantit\'e~$S$ se r\'eduit à une longueur rectiligne~$s$, il
suffira, pour le d\'emontrer dans le cas contraire, de d\'ecomposer
$S$ en \emph{\'el\'ements infiniment petits}.  (Cauchy [17], 1850,
p.\,171; emphasis added)
\end{enumerate}
As \emph{there is no conceptual analysis} of this issue and as their
putting ``infinitely small elements'' in italics proves, it is
sufficient for their intended appropriation of Cauchy as a forerunner
of non-standard analysis to find the term ``infinitesimal'' or
``infiniment petit'' in his texts.  \cite[p.\;4]{Sc22} (emphasis on
``there is no conceptual analysis'' added)
\end{enumerate}
Here GS claims that
\begin{enumerate}
\item
there is ``no conceptual analysis of this issue'' in
\cite{20b}, and that
\item
the conclusion is based merely on the use of the
term ``infinitely small elements'' by Cauchy.
\end{enumerate}
However, both claims are inexact.  Not only does \cite{20b} provide an
analysis to justify its conclusion, but there is an almost page-long
subsection~4.2 entitled ``Analysis of Cauchy's argument"
\cite[pp.\;132--133]{20b} where the conclusion is justified (along the
lines of the summary provided at the beginning of this section).

GS is free to disagree with the interpretation given in \cite{20b},
but he is less free to misrepresent it for the readers of MathSciNet.

\subsection{Infinitesimals and null sequences}

In a further comment on Bair et al.~\cite{20b}, GS claims that
\begin{enumerate}\item[]
[T]he paper reveals several points where the authors withdraw claims
made in earlier publications of the group'' \cite[p.\,1]{Sc22}
\end{enumerate}
He attempts to illustrate his claim by contrasting two passages:
\begin{enumerate}
\item
a passage from 2013 to the effect that ``In Cauchy, any variable
quantity~$q$ that does not tend to infinity is expected to decompose
as the sum of a given quantity~$c$ and an
infinitesimal~$\alpha$:~$q=c+\alpha$'' \cite[pp.\;900--991]{13a}, and
\item
a passage from 2020 to the effect that ``Cauchy's presentation of
infinitesimal techniques [in the calculation of the radius of
  curvature] contains no trace of the variable quantities or sequences
exploited in his textbooks in the definitions of infinitesimals''
\cite[p.\,135]{20b}.  
\end{enumerate}
Comparing the two passages, GS claims that
\begin{enumerate}\item[]
in [the 2013] presentation of the essentials of their claims, the
notion of null sequences was apodictically excluded.  In [2020],
\emph{the authors admit} such a notion for Cauchy, {\ldots}
\cite[p.\,1]{Sc22} (emphasis on ``the authors admit'' added)
\end{enumerate}
The problem with GS's claim is two-fold:
\begin{enumerate}
\item
Contrary to Schubring's claim, null sequences were indeed mentioned in
2013, in the following terms: ``Cauchy handles the said notion
[related to uniform convergence] using infinitesimals, including one
generated by the null sequence $\left(\frac1n\right)$''
\cite[p.\;891]{13a}.
\item
Null sequences were discussed in detail in 2011 (and elsewhere) in the
following terms: ``infinitesimals themselves are defined in terms of
variable quantities becoming arbitrarily small (which have often been
interpreted as null sequences).  Cauchy writes that such a null
sequence `becomes' an infinitesimal~$\alpha$\,'' \cite[p.\;428]{11b}
\end{enumerate}
By seeking to contrast the 2013 and the 2020 articles, Schubring
commits himself to evaluating the work of ``at least 22 mathematicians
and philosophers'' as a whole (see Section~\ref{s21}).  He attempts to
present the discussion of null sequences in 2020 as some kind of late
afterthought (``the authors admit, etc.''), but ends up
misrepresenting the work by the group.  His odd claim concerning ``the
authors withdraw[ing] claims made in earlier publications of the
group'' is groundless, and does not improve with repetition:
\begin{enumerate}\item[]
\emph{Admitting} that Cauchy had exploited the null-sequences notion
in his textbooks, the authors aim to verify whether their
infinitesimalist interpretation is at least viable in certain
applications, {\ldots} \cite[p.\;2]{Sc22} (emphasis on ``Admitting''
added)
\end{enumerate}
As documented above, Cauchy's notion of null sequence was a constant
presence in the work by the group, rather than some kind of late
``admission'' as claimed misleadingly by GS.

\subsection{Ambiguous continuity}
\label{s33b}

In the same vein, GS goes on to claim that
\begin{enumerate}\item[]
Regarding the concept of continuity, the paper under review shows
another and even more remarkable change of position and withdrawal.
{\ldots} They comment `that Cauchy's definition of continuity is, from
a modern viewpoint, somewhat ambiguous' (Bair et al.~2020, p.\;140).
All former certainty is gone. And this is \emph{not too far} from the
result of the careful analysis of Umberto Bottazzini, who chose to
speak of `ambiguous' and to attribute to Cauchy a special meaning of
continuity, `C-continuity' \cite[pp.\;4--5]{Sc22} (emphasis on ``not
too far'' added)
\end{enumerate}
Here GS juxtaposes the use of the adjective \emph{ambiguous} by Bair
et al.~\cite{20b} and by Bottazzini, and suggests that they are ``not
too far.''  The problem with GS's claim is two-fold:
\begin{enumerate}
\item
\label{bott}
Cauchy's final definition of continuity requires~$f(x+\alpha)$ to be
infinitely close to~$f(x)$ for all infinitesimal~$\alpha$.  The
definition is ambiguous since it is unclear whether Cauchy meant this
to apply at ``ordinary''~$x$ or, in addition, at~$x$ generated by
variable quantities.%
\footnote{For a modern formalisation see e.g., \cite{21d}.}
This is the ambiguity referred to by Bair et al.~in 2020 on page 140,
involving no change of position relative to earlier articles by the
group.
\item
Declaring Cauchy's continuity to be ``C-continuity'' (or for that
matter Weierstrass's continuity to be ``W-continuity'') carries as
little explanatory power as explaining the classical
infinitesimalists' work by their posession of ``an unerring
intuition'' \cite[p.\;358]{Gr74}.
\end{enumerate}
While this is not the place to analyze Bottazzini's concept, it is
necessary to point out that it is unrelated to the issue discussed in
item~\eqref{bott} above, contrary to GS's suggestion.  Thus GS's claim
of ``remarkable change of position and withdrawal'' is baseless.

\subsection{Their own non-standard-analysis concept}
\label{s86}

GS lodges the following claim concerning Cauchy's colleagues Poisson,
de Prony, and Petit:
\begin{enumerate}\item[]
[T]he authors affirm that all of them [i.e., Cauchy's colleagues] not
only shared the same conceptions, but even that there had been just
one unique conception of infinitesimals: \emph{their own
  non-standard-analysis concept}:
\begin{enumerate}\item[]
There seems to be little reason to doubt that the notion of infinitely
small in the minds of Poisson, de Prony, Petit, and others was solidly
in the Leibniz--l'H\^opital--Bernoulli--Euler school. (Bair et
al. 2020, p.\,142)
\end{enumerate}
\cite[p.\;6]{Sc22} (emphasis added)
\end{enumerate}
What exactly is the connection between GS's claim and the indented
quotation from Bair et al.~\cite{20b}?  The indented passage
mentioning Poisson, de Prony,%
\footnote{GS's stance on de Prony is puzzling.  Namely, GS attributes
  to de Prony ``the exclusion and rejection of \emph{infiniment
    petits} by the analytic method.  In de Prony the \emph{infiniment
    petits} were excluded from the foundational concepts of his
  teaching by simply not being mentioned; etc.''
  \cite[p.\;289]{Sc05}.  GS's claim flies in the face of de Prony's
  detailed treatment of the problem of infinitesimal oscillations, as
  well as his derivation of the formula for~$\cos z$ in terms of the
  exponential function following Euler; see \cite[Section~3.6]{19a}.}
Petit, Leibniz, l'H\^o\-pital, Bernoulli, and Euler%
\footnote{GS's stance on Euler is puzzling.  He claims that ``Euler
  established a purely algebraising foundation, achieving its climax
  in Lagrange's theory of functions'' \cite[p.\;1]{Sc22}.  The claim
  is meaningless without specifying what ``algebraising'' means
  exactly.  GS's intention here seems to be to minimize the importance
  of Euler's infinitesimals; related remarks were made by Ferraro (see
  item~\eqref{ferraro} in Section~\ref{s33}).  A rebuttal of such
  views apears in \cite{17b}.}
is apparently quoted as evidence of an alleged positing of a ``unique
{\ldots}\;non-standard-analysis concept.''  It does not require great
analytical skills on the part of the reader to ascertain that
nonstandard analysis was not mentioned at all in the passage quoted as
evidence for GS's non-standard claim.  Rather, the passage emphasized
the 17th--18th century infinitesimalist tradition.  GS's claim that
\cite{20b} attributed a ``non-standard-analysis concept'' to Cauchy's
colleagues is therefore as baseless as his earlier rhetorical
flourishes targeting Laugwitz quoted in Sections~\ref{s21} and
\ref{s212}.

\subsection{Hermetic monologue?}
\label{s212}

In 1989, Laugwitz presented a detailed comparative analysis of the
approaches of Fourier, Poisson, and Cauchy to the method of auxiliary
multipliers for obtaining values of indeterminate series and
integrals; see \cite[Section 5, pp.\;218--232]{La89}.  Thus, Laugwitz
was clearly aware of the fact that taking into account contemporary
work by Fourier and Poisson is essential for understanding Cauchy
himself.  One would not guess as much from the portrayal of Laugwitz'
work as depicted by GS:
\begin{enumerate}
\item
``[Laugwitz] reduced [Cauchy's universe of discourse] to a hermetic
  monologue by Cauchy" \cite[p.\;4]{Sc05};
\item
``they [Laugwitz and others] did not use Cauchy’s communication with
  contemporary mathematicians as a means to uncover what the
  respective concepts meant in their own period" \cite[p.\;4]{Sc05};
\item
``{\ldots} Laugwitz {\ldots}~had practically assigned a solipsistic
  mathematics to Cauchy" \cite[p.\;434]{Sc05};
\item
``{\ldots}~Laugwitz's approach of seeking meaning exclusively through
  the internal `conceiving' of a text is typical''
  \cite[p.\;441]{Sc05}.%
\footnote{While sharp criticisms of Laugwitz's Cauchy scholarship by
  Schubring and by Spalt are well known, what is perhaps less known is
  Schubring's sharp criticisms of Spalt's Cauchy scholarship; in fact
  all of the criticisms cited here target both Laugwitz and Spalt.}
\end{enumerate}
Given Laugwitz's attention to the historical context of Cauchy's work
as mentioned above, GS's ``hermetic monologue, etc.''  comments amount
to a strawman criticism.

Bair et al.~\cite{20b} commented in detail on the work of Cauchy's
contemporaries Poisson, de Prony, and Petit (see Section~\ref{s86}).
Related comments on Fourier had already appeared in the work
\cite[p.\;61]{13a} quoted by GS, and elsewhere.  One would not guess
as much from Schubring's strawman depiction of \cite{13a}:
\begin{enumerate}\item[]
[T]hey had understood `context' as referring just to other parts of
the same text, and were defending Detlef Laugwitz's earlier
interpretations of Cauchy, attributing to him a proper `universe',
independent and free of all relation with contemporary mathematics:
{\ldots} \cite[p.\;5]{Sc22}
\end{enumerate}
The article \cite{13a} indeed contained a defense of Laugwitz from
GS's strawman criticisms already in 2013, but neither Laugwitz nor any
of the scholars caricatured by GS in the above comments ever
attributed to Cauchy ``a proper `universe' independent and free of all
relation with contemporary mathematics.''

\subsection{Procedures, foundations, and misconceptions}
\label{s97}

GS's comments about Poisson, de Prony, Petit (see Section~\ref{s86})
and ``non-standard numbers'' (see Section~\ref{s23b}) reveal his
greatest misconception concerning formalisations of historical
mathematics.  He appears to believe that interpreting the work of
historical infinitesimalists in terms of Robinson's framework for
infinitesimal analysis necessarily amounts to attributing
``non-standard numbers'' to those 19th century authors.  If this meant
attributing some kind of ``anticipation'' of ultrafilters%
\footnote{For the benefit of the reader not familiar with the details
  of the construction of the hyperreal field~$\astr$, it may be useful
  to recall that~$\astr$ is obtained as a quotient~$\R^\N/\mathcal{U}$
  where~$\R^\N$ is the space of sequences of real numbers, and
  ~$\mathcal{U}$ is a nonprincipal ultrafilter on~$\N$.  For details
  see e.g., \cite{17f}.  Recently it turned out that ultrafilters are
  unnecessary for analysis with infinitesimals, which can be developed
  conservatively over ZF; see \cite{21e}.}
to 19th century authors, it would certainly amount to an absurd
interpretation.

However, formalisation of the work of historical infinitesimalists in
terms of modern theories of infinitesimals involves only providing
suitable proxies for their procedures and inferential moves (see e.g.,
\cite{20b}; on procedures see further in \cite{21a}).  Such
formalisations presuppose no relation to the foundational aspect of
the grounding of Robinson's infinitesimal analysis via a set-theoretic
construction (see also \cite{Cu88}).  The viability of applying NSA to
interpreting the procedures of the historical infinitesimalists
depends crucially on the procedure/foundation distinction.

The related procedure/ontology distinction was emphasized in our 2017
objection to GS's position, published five years ago
\cite[p.\,126]{17e}.  Five years later, GS's 2022 effort \cite{Sc22}
reveals little awareness of our objections.

Unlike Schubring and Spalt (see Section~\ref{spalt}), Detlef Laugwitz
clearly realized this distinction.  Thus, in his 1987 publication in
\emph{Historia Mathematica} \cite{La87} he carefully distinguished
between
\begin{enumerate}
\item
his analysis of Cauchy's procedures presented in Sections 1 through
14, and 
\item
his proposed models of Cauchyan infinitesimals in terms of
modern infinitesimal theories, in his Section 15.
\end{enumerate}
See further in \cite[Section~6.2]{17e}.%
\footnote{A failure to appreciate such a distinction led GS to accuse
  Laugwitz of attributing hyperreal numbers to Cauchy; see
  note~\ref{f4} and the main text there.}

The key insight is that Robinson's procedures provide better proxies
than Weierstrassian ones.  Thus, Robinson's standard/nonstandard
distinction%
\footnote{For example, numbers come in two varieties: standard and
  nonstandard.  Nonzero infinitesimals are necessarily nonstandard.}
is a proxy for Leibniz's assignable/inassignable distinction.  The
latter found 19th century echoes in the work of No\"el and Manilius
(see Section~\ref{s24}), ignored by GS.
%
%

Schubring is not the only historian insufficiently sensitive to the
dichotomy of procedures \emph{vs} foundations.  In his latest book,
L\"utzen writes:
\begin{enumerate}\item[]
In 1966 Abraham Robinson {\ldots} showed that it is possible to enrich
the real numbers by infinitesimals in a consistent way.  In the
resulting universe of non-standard analysis, one can apply the usual
rules of operation as with with real numbers except the Archimedean
property {\ldots} Robinson argued that his non-standard analysis
vindicated Leibniz', Euler's, and other earlier mathematicians'
calculations using infinitesimals.  This claim has been challenged, in
particular because Robinson's construction of his new non-standard
universe used modern methods that were far out of the reach of the
earlier mathematicians.  \cite[p.\,132]{Lu22}.
\end{enumerate}
Robinson's \emph{construction} indeed used modern methods, but his
\emph{procedures} nonetheless provide better proxies for the
inferential moves found in Leibniz, Euler, and others, just as
Robinson claimed.

\subsection{Interpretations, lores, and crusades}
\label{s211}

GS has pursued a different interpretation of Cauchy, and had already
expressed himself in his key publication concerning Laugwitz's alleged
attribution of hyperreal numbers to Cauchy (see Section~\ref{s21}),
and in 2016 concerning his opponents' alleged ``misconceptions''
\cite{Sc16}.  However, note the following five points.

\medskip
1.  GS arguably goes beyond a scholarly disagreement when he resorts
to accusing his opponents of allegedly ``leading a crusade against
historiography'' \cite[p.\,1]{Sc22} (see Section~\ref{f3b}) and
``falsifying Cauchy's text'' (see Section~\ref{s36}), as he did when
he accused Laugwitz of attributing hyperreal numbers to Cauchy (see
Section~\ref{s21}), as well as of attributing ``a hermetic monologue''
and ``solipsistic mathematics'' to Cauchy (see Section~\ref{s212}).

\medskip
2.  GS formulates a sweeping indictment not merely of the
article~\cite{20b} but of an entire research program (see
Section~\ref{s21}).  But then he goes on to accuse the authors of
\cite{20b} of failing to mention Cauchy's verb \emph{devient}, while
himself failing to mention the detailed analyses of \emph{devient} in
the earlier articles (see Section~\ref{s36}).  Such selective coverage
comes dangerously close to being mendacious by omission, and falls
short of a legitimate scholarly criticism.

\medskip
3.  GS's one-sided portrayal of 19th century \emph{infinitesimal
  lore}, where allegedly infinitesimals were necessarily ``standard
lore for expressing an arbitrarily small number'' \cite[p.\;3]{Sc22},
is contrary to historical fact; see the Belgian debate in
Section~\ref{s24} and the analysis there of the dual nature of 19th
century infinitesimal lore.  His one-sided portrayal is also contrary
to what GS wrote in his key publication concerning the
``far-reaching'' influence in France of Poisson's genuine
infinitesimals (or the ``pseudo-\emph{infiniment petits}'' in GS's
parlance), clearly establishing a plurality of such \emph{lores}.
Thus, what GS claims in 2022 is at odds with his own interpretation as
developed in his key publication.

\medskip
4.  GS's contention that Bair et al.~\cite{20b} attributed
``non-standard numbers'' to 19th century authors (see
Sections~\ref{s23b} and \ref{s86}) is a strawman criticism (parallel
to a similar contention with regard to Laugwitz mentioned in item 1).
Indeed, no such claim appeared in \cite{20b}, where it was emphasized,
on the contrary, that procedures based on modern infinitesimals are
only proxies for the procedures of the historical infinitesimalists
(see Section~\ref{s97}).  Furthermore, \cite{20b} clearly acknowledged
the profound differences in the background ontology of the historical
infinitesimalists and the modern ones.  There are further
misrepresentations in GS's 2022 review \cite{Sc22} but we will limit
ourselves to the remarks already made, for lack of space.

\medskip
5.  GS claims that when Cauchy's 1826 text on differential geometry
(see Section~\ref{s24b}) speaks of~$\varepsilon$ as ``un nombre
infiniment petit'' he is \emph{not} referring to a genuine
infinitesimal (or a \emph{pseudo-infiniment petit} in GS's parlance).
Here Cauchy develops a formula for the radius of curvature~$\rho$ at a
point of a curve in terms of the variation~$\Delta\tau$ of its tangent
vector.  Significantly, this particular~$\varepsilon$ occurs in the
same equation as Cauchy's increment~$\Delta\tau$:
\begin{equation}
\label{e61b}
\frac{\sin\left(\frac{\pi}{2}\pm\varepsilon\right)}{r}=
\frac{\sin(\pm\Delta\tau)}{\sqrt{\Delta x^2+\Delta y^2}}.
\end{equation}
From \eqref{e61b}, Cauchy derives the relation
\[
\frac{1}\rho=\pm \frac{d\tau}{\sqrt{dx^2+dy^2}}
\]
involving the radius of curvature~$\rho$, by passing to the limits
\cite[pp.\;98--89]{Ca26}.
%
%
Cauchy refers to his~$\Delta\tau$ as an \emph{angle de contingence}
(ibid.), sometimes translated as \emph{hornlike angle}.  This is a
traditional term for an angle incomparable with ordinary rectilinear
angles at least since the 16th century.  As GS acknowledged in his key
publication,
\begin{quote}
Klein has shown in detail that the hornlike angles form a model of
non-Archimedean quantities.  \cite[p.\,17]{Sc05}
\end{quote}
The reference is to \cite[p.\;221]{Kl25} (see also 
\cite{La73}).%
\footnote{On Klein see further in \cite{18b} and \cite{18i}.}
The occurrence of~$\varepsilon$ and the hornangle~$\Delta\tau$ in the
same equation indicates that this particular~$\varepsilon$ shares the
non-Archimedean nature of~$\Delta\tau$.  In more detail, since the
hornangle~$\Delta \tau$ on the right-hand side of \eqref{e61b} is
infinitesimal,~$\sin (\pm \Delta\tau)$ is also infinitesimal;
therefore~$\sin\left(\frac{\pi}{2}\pm\varepsilon\right)$ on the
left-hand side is infinitely close to~$1$; hence~$\varepsilon$ is also
infinitesimal.  In 2022, GS's account of Cauchy's 1826 text in
differential geometry fails to clarify the non-Archimedean context of
the~$\varepsilon$ as used by Cauchy.  Schubring's assumption that
$\Delta\tau$ is a finite difference betrays a presentist bias.

\subsection{Default number systems}
\label{s39}

Our gentle reader may well wonder why explorations of Cauchy's use of
genuine infinitesimals in work by Laugwitz and Bair et al.  should
provoke GS to describe such work as anachronistic attempts to depict
Cauchy as hermetic, solipsistic, and ``free of all relation with
contemporary mathematics'' (see Section~\ref{s212}).  This is
especially puzzling since GS himself clearly enunciated his
disagreement with Grabiner's portrayal of Cauchy as a pioneer of the
Weierstrassian \emph{Epsilontik} (see Section~\ref{f18}).
The explanation depends crucially on the distinction between the
following two questions:
\begin{enumerate}
\item[(A)] Did Cauchy's work contain a significant \emph{Epsilontik}
  component?
\item[(B)] Did Cauchy use only Archimedean quantities?
\end{enumerate}
The fact that GS identifies fully with 19th century critics of
infinitesimals such as Moigno and Mansion suggests that, while GS
would answer question~(A) in the negative, he assumes Cauchy's
background continuum to be Archimedean as the \emph{default} option
that requires no further argument.  Such an assumption (a natural
product of modern undergraduate training in naive set theory and
calculus/analysis) leads GS to assume further that Cauchy's
contemporaries would have been unable to understand Cauchy had the
latter used genuine infinitesimals.  Significantly, GS assigns
ontological import to Cauchy's expression \emph{absolute numbers},
whereas in reality Cauchy merely referred to the convention of
interpreting unsigned numbers~$x$ as positive ($+x$) rather than
negative ($-x$) numbers (see Section~\ref{f3b}).  Products of such
default thinking are the `hermetic' and `solipsistic' flourishes (see
Section~\ref{s212}); no wonder GS believes infinitesimals to be
``vaguely conceived'' (see Section~\ref{s13}).%
\footnote{In the same vein, Spalt appears to view genuine
  infinitesimals as bordering on the supernatural: ``an `infinitely
  small quantity' is for Leibniz nothing supernatural,
  inconceivable--but only a special case of a commonly used changing
  quantity: just one which decreases indefinitely''
  \cite[p.\;36]{Sp22}.}

Indeed, the outcome of Schubring's 2022 analysis in \cite{Sc22} is
predetermined by his historiographic assumptions, one of which seems
to prohibit Cauchy from using genuine infinitesimals, inspite of the
evidence provided by Laugwitz in his publications in \emph{Historia
  Mathematica} \cite{La87} and \emph{Archive for History of Exact
  Sciences} \cite{La89} and elsewhere, and by Bair et al.~in
\emph{British Journal for the History of Mathematics} \cite{20a} and
elsewhere.  We will analyze some related historiographic issues in
Section~\ref{s20}.

\section{Royal road to the great triumvirate}
\label{s20}

\subsection{Teleology through examples}
\label{s33}

One persistent theme in the historiography of the 19th century is the
perception that mathematical analysis reached its teleological
fulfillment with the development of what mathematicians often consider
to be ultimate foundations of analysis during the period of
Weierstrass and following.  Such perspectives typically include the
conception of an infinitesimal-free continuum as the true foundation
of analysis, and come assorted with an enduring faith in a literal
interpretation of the epithet \emph{real} in the expression \emph{real
  number} (an attitude that tends to overlook the fact that we only
have a \emph{theory} of real numbers, not an absolute standard
\emph{model} a.k.a.~the \emph{intended interpretation}),%
\footnote{Such attitudes are common among mathematicians who describe
  themselves as Platonists, such as Alain Connes; see further in
  \cite{13c}, \cite{13d} \cite{15a}, and \cite[Section~3.5]{18l}.}
and the accompanying enduring faith that the elimination of
infinitesimals was an inevitable part of the teleological process.
%
%
We provide some examples from the recent literature on Leibniz, Euler,
and Cauchy.

\begin{enumerate}
\item
There have been sustained efforts ranging from Ishiguro (\cite{Is90},
1990) to Rabouin and Arthur (\cite{Ra20}, 2020) to deny that Leibniz
meant the term \emph{infinitesimal} to refer to a mathematical entity
(see Section~\ref{s12}).  Most recently in 2022, Arthur persists in an
error already diagnosed by Fraenkel over a century ago; see
Section~\ref{s15}.%
%
%
\item
\label{ferraro}
Commenting on Euler's \emph{Introductio in analysin infinitorum},
Ferraro perceives a causal connection between the use of
infinitesimals and lack of success:
\begin{enumerate}
\item[] ``Euler was not entirely successful in achieving his aim since
  he introduced infinitesimal considerations in various proofs.''%
\footnote{As argued in \cite{17b}, Euler was more successful in
  achieving his aim than some historians believe.}
\cite[p.\,11]{Fe20}
\end{enumerate}
\item
Siegmund-Schultze suggests that infinitesimals in Cauchy are not
merely a remnant of the past but actually constitute a step backward:
\begin{enumerate}\item[]
``There has been \ldots{}\;an intense historical discussion in the
  last four decades or so how to interpret certain apparent remnants
  of the past or -- as compared to J. L. Lagrange's (1736--1813)
  rigorous `Algebraic Analysis' -- even \emph{steps backwards} in
  Cauchy's book, particularly his use of infinitesimals {\ldots}''
  (\cite{Si09}; emphasis added)
\end{enumerate}
\end{enumerate}

Attempts to present Cauchy as a precursor of Weierstrass in the
received Cauchy literature are based on similar assumptions,%
\footnote{See Section~\ref{s21b}.}
and ignore numerous works and applications where Cauchy used
infinitesimals as numbers, as detailed in~\cite{20b} (see
Section~\ref{s24} for an example).%

Schubring claims that the group Bair et al.
\begin{enumerate}\item[]
has been trying to rewrite the history of the infinitesimal calculus
as a forerunner of non-standard analysis \cite[p.\,1]{Sc22}
\end{enumerate}
but overlooks the fact that Bair et al.~only object to
\begin{enumerate}
\item
writing the history of the infinitesimal calculus as necessarily a
forerunner of the Weierstrassian \emph{Epsilontik}, and
\item
Schubring's assumption that Cauchyan infinitesimals belong in an
infinitesimal lore (see Section~\ref{s24}) that is unequivocally
Archi\-me\-dean.
\end{enumerate}
Schubring's ``old-fashioned'' \cite[p.\,1]{Sc22} take on Leibniz is
analyzed in Section~\ref{s13}, and his take on Cauchy,
in~Sections~\ref{s2} and \ref{s3}.

\subsection{Fraser on inevitable evolution}

Certainly, some historians are aware of the pitfalls of teleological
fallacies analyzed in Section~\ref{s33}.  Thus, Fraser writes:
\begin{enumerate}\item[]
Since the 1960s there has been a new wave of writing about the history
of eighteenth-century mathematics.  Authors such as Henk Bos, Steven
Engelsman, Niels Jahnke, Giovanni Ferraro, Craig Fraser and Marco
Panza have charted the development of calculus without interpreting
this development as a first stage in the inevitable evolution of an
arithmetic foundation.  \cite[p.\;27]{Fr15}
\end{enumerate}
Fraser appears to acknowledge that the traditional ``arithmetic
foundation'' of classical analysis as developed around 1870 was not
the inevitable result of the evolution of analysis.  Yet he
immediately goes on to reassure his readers in the following terms:
\begin{enumerate}
\item
Of course, classical analysis developed out of the older subject and
it remains a primary point of reference for understanding the
eighteenth-century theories.  \cite[p.\;27]{Fr15} 
\item
The relevance of modern non-Archimedean theories to an historical
appreciation of the early calculus is a moot point.
\cite[p.\;43]{Fr15}
\end{enumerate}
Postulating the superiority of classical analysis over non-Archimedean
theories as the basis for a historical appreciation of the early
calculus involves precisely the type of teleological fallacy examined
in Section~\ref{s33}.  Fraser proceeds to succumb to a closely related
fallacy of conflating procedures and foundations (see
Section~\ref{s97}):
\begin{enumerate}\item[]
[N]onstandard analysis and other non-Archimedean versions of calculus
emerged only fairly recently in somewhat \emph{abstruse mathematical
  settings} that bear little connection to the historical developments
one and a half, two or three centuries earlier.  \cite[p.\;27]{Fr15}
(emphasis added)
\end{enumerate}
For the benefit of the reader not familiar with foundational
subtleties, we hasten to point out that Robinson's framework is
grounded in the traditional Zermelo--Fraenkel set theory, no longer
considered abstruse by classically-trained mathematicians.  We can
agree with Fraser that the \emph{foundational} aspects of grounding
Robinson's infinitesimals in classical set theory bear little
connection to the historical developments from Leibniz to Cauchy.
However, the \emph{procedures} of these pioneers of analysis do
exhibit a strong connection to those developed by Robinson.%
\footnote{See further in \cite[Sections~4.2--4.6, pp.\,123--128]{17a}
  and \cite[Section~4.4, pp.\;277--278]{18e} for an analysis of
  Fraser's text.}  Schubring and L\"utzen, as well, as insufficiently
sensitive to this distinction; see Section~\ref{s97}.

A recent piece by Archibald et al.~published in \emph{The Mathematical
  Intelligencer} in response to the article ``Two-track depictions of
Leibniz's fictions'' \cite{22b} enables a considerable extension of
the list of scholars who apparently have difficulty separating the
contention that Leibniz and others exploited procedures using genuine
infinitesimals, from the idea of a ``pervasive presence of nonstandard
analysis in the history of mathematics.''  The overlap between
Archibald's coauthors and Fraser's list of leading scholars (see
above) includes Ferraro and Panza, apparently providing evidence
against Fraser's claim quoted at the beginning of this section.  See
Section~\ref{s1v} for more details.

\subsection{Grattan-Guinness on presentism}
\label{s31}

In an influential 1990 article, Grattan-Guinness wrote:
\begin{enumerate}\item[]
[Mathematicians] usually view history as a record of a `royal road to
me' -- that is, an account of how a particular modern theory arose out
of older theories instead of an account of those older theories in
their own right.  In other words, they confound the question, `How did
we get here?', with a different question, `What happened in the past?'
\cite[p.\,157]{Gr90b} (cf.\;\cite[p.\,165]{Gr04})
\end{enumerate}
The critique of the ``How did we get here?'' attitude is on target.%
\footnote{Meanwhile, Grattan-Guinness's stereotyping of mathematicians
  is unfortunate as it comes dangerously close to endorsing Unguru
  polarity; see further in \cite{20e}.}
%


Grattan-Guinness's \emph{royal road to me} issue seems closely related
to the issue of presentism.  As an example, Grattan-Guinness
criticized the approaches of Dieudonn\'e and Birkhoff for presentism,
but painted a more sympathetic picture of the approach of Andr\'e
Weil:%
\footnote{Subsequently Grattan-Guinness was more critical of Weil in
  \cite[p.\,166]{Gr04}.}
\begin{enumerate}\item[]
For better sensitivity to the issues exhibited by another eminent
mathematician, see the lecture delivered by A.\,Weil in 1978 to an
international congress of mathematicians: ``History of Mathematics -
why and how'', {\ldots} \cite[note 19, p.\,170]{Gr90b}
\end{enumerate}



Grattan-Guinness's criticism does apply to Bourbaki's approach (to
writing the history of mathematics), and specifically in reference to
its teleological aspect; see further in Section~\ref{s33}.

\subsection{Hacking's dichotomy}
\label{s32}

In his last book, Ian Hacking (\cite{Ha14}, 2014) presented an
analysis closely related to Grattan-Guinness' critique.  Hacking
introduces a distinction between the butterfly model and the Latin
model for the development of a scientific discipline.  Hacking
contrasts a model of a deterministic (genetically determined)
biological development of animals like butterflies (the
egg--larva--cocoon--butterfly sequence), with a model of a contingent
historical evolution of languages like Latin.  Hacking notes that
\begin{enumerate}\item[]
If analysis had stuck to infinitesimals in the face of philosophical
naysayers like Bishop Berkeley, analysis might have looked very
different.  Problems that were pressing late in the nineteenth
century, and which moved Cantor and his colleagues, might have
received a different emphasis, {\ldots}\;%
%
%
This alternative mathematics might have seemed just as `successful',
just as `rich', to its inventors as ours does to us.
\cite[p.\,119]{Ha14}.
\end{enumerate}
To borrow Hacking's terminology, one could say that some historians of
mathematics seem convinced that the butterfly of rigorous analysis
needed to shed its infinitesimal cocoon in order to fly.  Emphasizing
determinism over contingency can easily lead to anachronism; see
further in~\cite{19a}.  Dauben asks:
\begin{enumerate}\item[]
Is it anachronistic to use nonstandard analysis or transfinite numbers
to ``rehabilitate'' or explain the works of Leibniz, Euler, Cauchy, or
Peirce, for example, as recent mathematicians, historians, and
philosophers of mathematics have attempted?  \cite[p.\;307]{Da21}
\end{enumerate}
He answers as follows:
\begin{enumerate}\item[]
Robinson succeeded in showing the reasonableness of ``redrawing'' the
early history of the calculus to reinstate past views that, cast in
the light of nonstandard analysis, could be seen more clearly.  In
these cases at least, an anachronistic explanation nevertheless serves
to clarify, not confound, what had confused earlier defenders of
theories based on infinitesimals like the calculus.
\cite[p.\;327]{Da21}
\end{enumerate}
Dauben does not analyze the issue of anachronism in terms of the
procedure \emph{vs} foundation distinction (see Section~\ref{s97}),
but arguably an explanation addressing the procedures while
acknowledging the differences in foundations, no longer needs to be
described as anachronistic if it succeeds in accounting for the
inferential moves of the historical authors.

\subsection{Contingency and determinism}

The contingency of the historical evolution of the mathematical
sciences would entail in particular that the mathematical landscape
today could have been different from what it currently is.  Such a
perspective is consonant both with Hacking's Latin model (see
Section~\ref{s32}) and with Grattan-Guinness's critique of the ``How
did we get here?''  approach (see Section~\ref{s31}).  The shortcoming
of the latter approach is its implied faith in the determinism of the
historical evolution of mathematics.
In the history of any science, it is philosophically problematic to
claim a singular juncture when suddenly what was obscure becomes
clarified.%
\footnote{A telling example of a postulation of such a singular
  juncture, following a millenial aspiration, is provided by Frank
  Quinn: ``The breakthrough was development of a system of rules and
  procedures that really worked, in the sense that, if they are
  followed very carefully, then arguments without rule violations give
  completely reliable conclusions.  {\ldots} There is no abstract
  reason (i.e., apparently no proof) that a useful such system of
  rules exist, [sic] and no assurance that we would find it.  However,
  it does exist and, after \emph{thousands of years} of tinkering and
  under intense pressure from the sciences for substantial progress,
  we did find it'' \cite[pp.\;31--32]{Qu12} (emphasis on ``thousands
  of years'' added).  To illustrate the touted millenial breakthrough,
  Quinn provides the example of ``Weierstrass's nowhere-differentiable
  function (1872)'' \cite[p.\;31]{Qu12}.  We will not comment on
  Quinn's assumption that the rigorisation of mathematical analysis
  occurred under ``intense pressure from the sciences.''}

Some historians of mathematics tend to view the history of
mathematical analysis as an exception in this regard, possibly under
the influence of their undergraduate training in (naive) set theory
and the \emph{Epsilontik}.%
\footnote{\label{f37}Thus, Knobloch \cite[pp.\;13--14]{Kn18} does not
  hesitate to appeal to Alef$_0$ [his notation] in a discussion of the
  Leibnizian calculus.  See further in
%
\cite[Section~3]{21g} and \cite[Section~3.3]{23c}.}
Such historians pursue a teleological reading of the history of
analysis of the 17--19 centuries, giving credit to the ``great
triumvirate'' \cite[p.\;298]{Bo49} of Cantor, Dedekind, and
Weierstrass in this connection, in a version of the approach mocked by
Grattan-Guinness that can be characterized as a \emph{royal road to
  the great triumvirate}.  In general, historians of natural science
make no analogous ``singular juncture'' claims, and are on more solid
ground historiographically speaking.  As noted by Gray, ``in
mathematics, as in the rest of science, authority is only partial,
dynamic, and contested'' \cite[p.\;512]{Gr11}.%
\footnote{In the same text \cite[p.\;514]{Gr11}, Gray unfortunately
  also endorses Mehrtens' odd compendium of misinformation \cite{Me90}
  on Felix Klein (as Gray already did in his book \cite{Gr08a}).  We
  set the record straight in \cite{18b}.}

%

Viewed from this perspective, it is perhaps illusory, or more
precisely circular, to insist that concepts like \emph{infinitesimal}
must be clarified before they can serve a useful purpose in an
argument.  Both the concepts and the arguments are part of an evolving
effort by scholars to reach clearer understanding.  Arguably this
contention applies as much in mathematics as it does in natural
science; the singular emphasis placed in some history of mathematics
books on the foundational developments of the 1870s may amount to such
a \emph{royal road}.  For all the seminal \emph{mathematical}
importance of the developments of the 1870s (known to all), their
\emph{philosophical} force was undercut by the mathematical
developments of the 1970s when Karel Hrbacek \cite{Hr78} and Edward
Nelson \cite{Ne77} developed axiomatic approaches to analysis with
infinitesimals; see \cite{23c} for further details.
%
%

\section{Is pluralism in the history of mathematics possible?}
\label{s1v}

As elaborated in Sections~\ref{s1} through \ref{s3}, the authors of
the present article have over the years developed the following
perspective.  Many mathematicians in the 17--19th centuries (Newton,
Leibniz, Euler, Cauchy, \ldots) employed one version or another of
informal calculus with infinitesimals.  When adhering to the internal
rules of such calculi, they produced valid results and predictions in
mathematics and physics.


Abraham Robinson's Nonstandard Analysis (NSA for short; 1961),
building upon earlier work by T. Skolem \cite{Sk33}, E. Hewitt
\cite{He48}, J. {\L}o\'s \cite{Lo55} and others, was a milestone in
that it constituted the first fully rigorous approach to `calculus
with infinitesimals'.  In light of the groundbreaking nature of
Robinson's work, the role of pre-Robinson/pre-20th century
infinitesimal calculi in the history of mathematics could be viewed as
follows.
\begin{enumerate}
\item
The procedures of pre-20th century infinitesimal calculi are
formalized more successfully in NSA than in traditional mathematics
based on the epsilon-delta framework in a purely Archime\-dean
context.
\item
The pre-20th century infinitesimal calculi should (at most) be viewed
as an inspiration for NSA.\, In particular, classifying them as
`forerunners' of NSA (or similar concepts) is problematic in that such
a classification projects a modern viewpoint onto 17--19th century
mathematics.
\end{enumerate}

We stand by this position, and in this section we defend the viability
of using NSA in this sense to interpret the procedures of the
historical infinitesimalists, against a recent attack by Tom Archibald
et al.~\cite{Ar22b}.  Furthermore, we feel that the possibility of
pluralism in the historiography of mathematics has just received an
unwarranted blow from Archibald et al.



\subsection{Depictions put in the pillory}

The article ``Two-track depictions of Leibniz's fictions'' \cite{22b}
was published in the september issue of the 2022 volume of \emph{The
  Mathematical Intelligencer}.  ``Two-track depictions'' analyzed
rival interpretations of the procedures of the Leibnizian calculus,
one of the issues being whether or not Leibniz used genuine
infinitesimals.  As emphasized in an earlier article \cite{21a} in the
\emph{British Journal for the History of Mathematics} and elsewhere,
an analysis of the procedures of the historical infinitesimalists
needs to be carefully distinguished from the foundational issues of
the grounding of infinitesimals in modern set-theoretic frameworks
(see Section~\ref{s97}).

A response by Archibald et al.~to ``Two-track depictions'' was
published online in \emph{The Mathematical Intelligencer} on 29
september 2022; see \cite{Ar22b}.  Our brief reply appeared at
\cite{23a}.

The piece by Archibald et al.~reveals that some historians of
mathematics apparently have difficulty separating the contention that
Leibniz and others exploited procedures using genuine infinitesimals,
from the idea of a ``pervasive presence of nonstandard analysis in the
history of mathematics.''  Thus, Archibald et al.~claim the following:
\begin{enumerate}\item[]
The aim [of ``Two-track depictions'' and other articles] is
{\ldots}~to put various scholars in the pillory -- with accompanying
abusive epithets -- for not enthusiastically recognizing the pervasive
presence of nonstandard analysis in the history of
mathematics.~\cite[p.\,1]{Ar22b}
\end{enumerate}
Indeed, there was no presence (pervasive or otherwise) of NSA in the
history of mathematics before 1961 when it was first introduced by
Robinson in \cite{Ro61}, but the \emph{procedures} of NSA provide
better proxies for the procedures of the historical infinitesimalists
than the procedures of Weierstrassian analysis in a purely Archimedean
setting; see further in~\cite{17d}.

Archibald et al.~do not provide any examples of alleged ``abusive
epithets'' in ``Two-track depictions'' -- for the simple reason that
there are none -- but see Section~\ref{ep}.

\subsection{Epithets and crusades}
\label{ep}

The piece by Archibald et al.~was signed among others by Gert
Schubring, the nature of whose epithets can be gleaned from comments
he made about the work of Cauchy historian Detlef Laugwitz, sampled in
Section~\ref{s212}.  Given Laugwitz's attention to the historical
context of Cauchy's work in his articles in \emph{Historia
  Mathematica} \cite{La87}, \emph{Archive for History of Exact
  Sciences} \cite{La89}, and elsewhere, Schubring's comments about
``hermetic monologue'' and ``solipsistic mathematics'' amount to a
strawman criticism already analyzed in \cite{17e}.

It is no secret that, like Laugwitz, the authors of \cite{21a} and
\cite{22b} have pursued an interpretation of Cauchy at variance with
Schubring's, in such venues as \emph{Perspectives on Science}
\cite{11b}, \emph{British Journal for the History of Mathematics}
\cite{20a} and elsewhere.  Schubring's reaction to such work is on
record.  Referring to the work of ``at least 22 mathematicians and
philosophers'' \cite{Sc22},%
\footnote{Schubring's unrefereed opinion piece is endorsed by
  Archibald et al.~\cite[note\;6]{Ar22b}.}
Schubring used the epithet ``crusade against the historiography'' to
describe such work (see Section~\ref{s21}).  Given such language, it
is decidedly comical to find Schubring accusing scholars -- who happen
to disagree with his historical interpretations -- of allegedly using
``abusive epithets.''

\subsection{Sole aim?}

Our aims in interpreting historical infinitesimalists were outlined at
the beginning of Section~\ref{s1v}.  Archibald et al.~claim that
\begin{enumerate}\item[]
these criticisms [contained in ``Two-track depictions'' and other
  articles] do not consist in reasoned historical argument, but rather
in piecemeal confutation of isolated quotations of their opponents
taken out of context, whose \emph{sole aim} is to show them as enemies
of the group's understanding of infinitesimals.  \cite[p.\,1]{Ar22b}
(emphasis added)
\end{enumerate}
Contrary to such claims, the ``22 mathematicians and philosophers''
(see Section~\ref{ep}) did present reasoned historical arguments
against what they see as untenable received interpretations, and
proposed better alternatives.

Consider for example the case of Richard Arthur, one of Archibald's
coauthors.  The 2021 article in \emph{British Journal for the History
  of Mathematics} \cite{21a} argued in detail that
\begin{enumerate}
\item
Arthur's attempt in \cite{Ar13} to interpret the Leibnizian calculus
in terms of a modern theory of infinitesimals called \emph{Smooth
  Infinitesimal Analysis} is unviable, and
\item
proposed a more viable alternative (which is currently the target
of~\cite{AR}).
\end{enumerate}

Another coauthor of Archibald's, Jeremy Gray, wrote the following
about Euler's foundations:
\begin{enumerate}\item[]
Euler's attempts at explaining the foundations of calculus in terms of
differentials, which are and are not zero, are \emph{dreadfully weak}.
\cite[p.\;6]{Gr08b} (emphasis added)
\end{enumerate}
Unfortunately, Gray provided no context for his ``dreadfully weak''
claim concerning Euler's foundations, suggesting that he assumes such
views to be generally accepted.  Gray's assumptions were challenged in
a detailed study of Euler in \emph{Journal for General Philosophy of
  Science}~\cite{17b}.  Since there has been no follow-up by Gray, it
is difficult to know how he would defend his assumptions concerning
Euler.

The Euler scholarship of Archibald's coauthor Giovanni Ferraro was
analyzed in detail in \cite{17b}, as well, similarly without
follow-up, except perhaps for the following recent comment.
Commenting on Euler's \emph{Introductio in analysin infinitorum},
Ferraro perceives a causal connection between the use of
infinitesimals and lack of success:
\begin{enumerate}\item[]
Euler was not entirely successful in achieving his aim since he
introduced infinitesimal considerations in various proofs.
\cite[p.\,11]{Fe20}
\end{enumerate}
Ferraro's 2020 comment is disappointingly consistent with the
presentist attitude in his earlier work analyzed in \cite{17b}.

\subsection{Gray and L\"utzen on Cauchy}

Gray also claimed the following concerning Cauchy's definitions:
\begin{enumerate}\item[]
[Cauchy] defined what it is for a function to be integrable, to be
\emph{continuous}, and to be differentiable, using careful, if not
altogether unambiguous, limiting arguments.  \cite[p.\;62]{Gr08a}
(emphasis added)
\end{enumerate}
Such a claim is inaccurate at least with regard to continuity, a point
argued in \cite{17b}.  A more successful alternative has been
elaborated in a number of publications on Cauchy including a 2020
article in \emph{British Journal for the History of
  Mathematics}~\cite{20a}.

Archibald's coauthor Jesper L\"utzen has been more forthcoming than
Gray with information about Cauchy.  Note the following five points
(summarizing the analysis in \cite[Section 3]{17a}):

\begin{enumerate}
\item
L\"utzen acknowledges that Cauchy's definitions contain no
quantifiers, writing: ``We miss our quantifiers, our $\varepsilon$'s,
$\delta$'s'' \cite[p.\;161]{Lu03}.
\item
He acknowledges that Cauchy's second definition of continuity used
infinitesimals \cite[p.\;160]{Lu03}.
\item
However, L\"utzen misrepresents the work of Robinson and Laugwitz when
he claims that they asserted that Cauchy's variables go through
infinitesimal values on their way to zero \cite[p.\;164]{Lu03}.
Neither Robinson nor Laugwitz ever made such a claim to our knowledge,
though it is found in a 1978 paper on Cauchy by Fisher
\cite[p.\;316]{Fi78}.
\item
He claims that the truth is found in Grabiner \cite{Gr81}, who
explains that whatever L\"utzen and others ``miss'' (see item (1)) is
actually found in Cauchy's \emph{proofs} (rather than definitions),
which are ``strikingly modern'' \cite[p.\;161]{Lu03}.
\item
According to L\"utzen, it is a ``fundamental lacuna'' of Cauchy's
proof of intermediate value theorem that the result relies on
completeness, which could not have been provided by Cauchy
\cite[pp.\,167--168]{Lu03}.  But as Laugwitz already pointed out
\cite[p.\;202]{La89}, Cauchy did not need a construction of the reals
because he had unending decimal expansions (available ever since Simon
Stevin).  Criticizing Cauchy's proof on the grounds of the missing
property of completeness therefore risks being anachronistic.
\end{enumerate}

A significant point concerns the disagreement between a pair of
Archibald's coauthors: L\"utzen endorses Grabiner's analysis of
Cauchy, whereas Schubring criticizes Grabiner's approach, as detailed
in Section~\ref{f18}.%
\footnote{Schubring's position is closer to Grattan-Guinness's, who,
  decades earlier, ``warn[ed] against planting later refinements
  (especially the Weierstrassians') onto that period [of Cauchy's
    activity in the 1820s]'' \cite[note 1, p.\;713]{Gr90}.}
The L\"utzen--Schubring disagreement leads us to an interesting
question of who exactly is entitled to disagree without running the
risk of being branded a crusader (see Section~\ref{ep}).

The list could be continued, but we hope to have illustrated the fact
that believing the sweeping claims by Archibald et al.~would entail
canceling large parts of modern scholarship published in leading
history and philosophy journals.

\subsection{Leibniz--Bernoulli correspondence}

Archibald et al.~make a number of additional spurious claims,
including the claim of having detected a contradiction in the work of
the ``22 mathematicians and philosophers'' (see Section~\ref{ep}).
Since Archibald's coauthors Arthur and Rabouin elaborate in
\cite{Ra20} on a spurious claim of having detected contradictions also
in the notion of infinitesimal in Leibniz (as analyzed in ``Two-track
depictions'' and in more detail in~\cite{23d}),%
\footnote{The alleged contradiction results from their tendency to
  assimilate infinitesimals to infinite wholes: essentially Rabouin
  and Arthur are trying to invert a cardinality to obtain an
  infinitesimal, a procedure that shocked Fraenkel over a century ago;
  see Section~\ref{s15}.}
we are satisfied to be in good company.

For the reader interested in the technical details, note that
Archibald et al.~claim to have detected a contradiction in our
interpretation of Leibniz.  The claim is based on Archibald et al.'s
reading of the Leibniz--Bernoulli correspondence from 1698--1699:
\begin{enumerate}\item[]
In the case [Leibniz] discussed with Bernoulli in 1698--1699, the
question was rather about whether there exists an infinitieth term in
an infinite series, which in the case of a decreasing series would
stand for an infinitesimal quantity.  Bernoulli insisted that there
would indeed be such an infinitieth (although not necessarily last)
term, thus entailing the existence of an infinitesimal. According to
Katz et al., this entails a conception of infinite series as
consisting of an infinite sequence of standard numbers followed by an
infinitesimal part.  \cite[p.\;2]{Ar22b}
\end{enumerate}
However, the claim by Archibald et al.~is based on a misreading of
crucial aspects of the 1698--1699 correspondence.  Basically,
Archibald et al.~are committing an elementary logical error, as we
will now explain.

In their correspondence, Bernoulli tried to convince Leibniz, through
the analysis of the behavior of an infinite series, that an
infinitesimal term in the series must exist.  To put it another way,
Bernoulli tried to derive the existence of infinitesimals from the
existence of infinite series.  Archibald et al.~appear to believe
that, since we referred to a non-Archimedean continuum as a
\emph{Bernoullian continuum}, we must agree with Bernoulli's
reasoning.  But note the following two points:
\begin{enumerate}
\item
The fact that Leibniz did not agree with Bernoulli's reasoning does
not mean that Leibniz rejected infinitesimals as (fictional)
mathematical entities; he merely found Bernoulli's argument flawed
because it was based on a conflation of magnitude and multitude.%
\footnote{A related misunderstanding occurs in Rabouin and Arthur
  \cite{Ra20}; for details see \cite{23d}, note 24, pp.\,12--13 and
  the main text there.}
\item
We similarly don't agree with Bernoulli's reasoning, and used the term
\emph{Bernoullian continuum} only because Bernoulli routinely used
infinitesimals in his \emph{mathematical} work (and not because we
agree with his reasoning from series in favor of the existence of
infinitesimals), quite apart from his \emph{philosophical} attempts to
convince Leibniz to adopt a more realistic position with regard to
infinitesimals.
\end{enumerate}
We have addressed this point concerning Bernoulli in detail because
Archibald et al.~apparently attach great importance to it, seeing that
the name ``Bernoulli'' is mentioned no fewer than 15 times in their
4-page text.

\subsection{Alice and Bobs}

With regard to the Leibnizian calculus, the 2022 article ``Two-track
depictions'' \cite{22b} presented and compared two interpretations,
represented respectively by Alice and Bob.  One of the aspects of
Leibniz's position highlighted by Bob was the notion of
\emph{infinitum terminatum} (lit.~bounded infinity), contrasted by
Leibniz with \emph{infinitum interminatum} (unbounded infinity).

Leibniz's position, as explained by the Leibniz historian Eberhard
Knobloch \cite[p.\;97]{Kn99}, is that the \emph{infinitum
  interminatum}, corresponding to an infinite whole (such as an
unbounded infinite line) is a contradictory notion (see
Section~\ref{s13}), whereas, by contrast, the \emph{infinitum
  terminatum} is a notion useful in geometry and calculus.  Such a
bounded infinity, as the name suggests, is exemplified by a subline
(bounded by a pair of infinitely separated endpoints) of the
(contradictory) unbounded infinite line.

Archibald et al.~quote this 2022 article, as well as the 2012 article
\cite{14c} and the 2013 article \cite{13f}.  Such attention is surely
appreciated by every researcher; hopefully it can mark the beginning
of a meaningful dialog or informed debate.  Archibald et al.~proceed
to label the three Bob2012, Bob2013, and Bob2021 (the latter seems to
be a misdated reference to the 2022 article), and to claim that
Bob2012 and Bob2021 contradict each other in their opinion of whether
infinitesimals are contradictory notions or not.  However, the
formulation ``contain a contradiction" in \cite{14c} (concerning
infinitesimals but also negatives and imaginaries) was in a different
context and must not be conflated with the contradictory nature of
Leibnizian ``infinite wholes.''  Toward the end of \cite{Ar22b}, one
finds an interesting footnote 4 to the effect that
\begin{enumerate}\item[]
Richard Arthur and David Rabouin, two of the authors of this paper,
will dedicate a specific study to this, providing several sources in
which Leibniz explicitly claimed that \emph{lineae infinitae
  terminatae} are contradictory entities.  \cite[note 4]{Ar22b}
\end{enumerate}
At the very least, it seems that Arthur and Rabouin owe thanks to Bob
for raising such an interesting issue, if it led to a new ``specific
study'' of theirs.  We look forward to seeing their ``specific study''
and suspend judgment of the merits of the, frankly surprising, claim
that the \emph{infinitum terminatum} is a contradictory notion -- a
claim that does not square with several texts by Leibniz where the
usefulness of the \emph{infinita terminata} is contrasted with the
contradictory nature of the \emph{infinita interminata}, as documented
in recent articles \cite{21a} and \cite{23d}.  Surprisingly, Archibald
et al.~expect the reader to accept a wholesale dismissal of the
research of ``22 mathematicians and philosophers'' (to quote
Schubring), published in leading journals, based on an
as-yet-unpublished specific study by Arthur and Rabouin.

\subsection{Pacidius}

There is a Leibnizian passage from 1676 where both \emph{infinita
  terminata} and contradictions are mentioned.  However, the passage
leads to the opposite conclusion from the one sought by Arthur and
Rabouin:
\begin{enumerate}\item[]
``Pacidius: I would indeed admit these infinitely small spaces and
  times in geometry, for the sake of invention, even if they are
  imaginary.  But I am not sure whether they can be admitted in
  nature.  For there seem to arise from them infinite straight lines
  bounded at both ends, as I will show at another time; which is
  absurd.''%
\footnote{``Ego spatia haec et tempora infinite parva in Geometria
  quidem admitterem, inventionis causa, licet essent imaginaria.  Sed
  an possint admitti in natura delibero. Videntur enim inde oriri
  lineae rectae infinitae utrinque terminatae, ut alias ostendam; quod
  absurdum est'' \cite[p.\;206]{Ar01}.}
(Leibniz as translated by Arthur in \cite[p.\;207]{Ar01})
\end{enumerate}
The structure of Leibniz's argument, consistent with his fictionalist
views developed elsewhere, is that there are no infinitesimals in
nature because if there were some, then there would also be
\emph{infinita terminata}, which would be absurd (by an argument that
Leibniz promises to provide elsewhere).  Accordingly, it is the
hypothesis of the existence of bounded infinities \emph{in nature}
that leads to an absurdity.  On the other hand, their usefulness in
geometry does not depend on their existence in nature.  This is a
powerful argument against the Rabouin--Arthur interpretation.

More generally, it is surprising that Archibald et al.~should present
the position of Arthur and Rabouin concerning Leibniz as allegedly
universally accepted among Leibniz scholars.%
\footnote{\label{f5}Thus, Archibald et al.~claim that ``there is no
  B-methodology sensu stricto in Leibniz. Leibniz's main argument is
  that it is not possible to treat infinitesimals as existing entities
  because that amounts to the introduction of an infinite number,
  which he takes to be a contradictory notion'' and go on to describe
  such a position as ``fact'' in their note 3: ``$^3$This had been a
  well-known fact among Leibniz scholars for some time''
  \cite[note\;3]{Ar22b}.  One of the works they cite is Bassler
  \cite{Ba98}.  The shortcomings of Bassler's reading are analyzed in
  \cite[Section\;2]{21g}.}
Quite the contrary: the 2020 article by Rabouin and Arthur on Leibniz
in \emph{Archive for History of Exact Sciences} was followed in the
same journal by the 2021 article by Esquisabel and Raffo Quintana
\cite{Es21}, who explicitly reject the Rabouin--Arthur interpretation
in the following terms:
\begin{enumerate}
\item
``[U]nlike the infinite number or the number of all numbers, for
  Leibniz infinitary concepts do not imply any contradiction, although
  they may imply paradoxical consequences.''  \cite[p.\;641]{Es21}
\item
``[W]e disagree with the reasons [Rabouin and Arthur] gave for the
  Leibnizian rejection of the existence of infinitesimals, and in our
  opinion the texts they refer to in order to support their
  interpretation are not convincing.  Since we argue that Leibniz did
  not consider the concept of infinitesimal as
  \emph{self-contradictory}, we try to provide an alternative
  conception of impossibility.''\;\cite[p.\;620]{Es21}
\end{enumerate}
Curiously, the article by Esquisabel and Raffo Quintana was
communicated by no other than{\ldots} Archibald's coauthor Jeremy Gray
(see \cite[p.\;613]{Es21}).  Apparently, Gray did not read carefully
one of the two texts: either the article he communicated, or the piece
by Archibald et al.~before consenting to have his name added to its
author list.


\subsection{Two methods in Leibniz}
\label{jesseph}

Perhaps the most remarkable case of membership on the author list of
Archibald et al.~is Douglas Jesseph.  Speaking of the law of
continuity in 1989, Jesseph asserts that
\begin{enumerate}\item[]
Leibniz argues that, when applied to the calculus, this law yields a
\emph{new kind of quantity} which will provide the foundation for the
reasonings which appear in the solution to geometrical problems.
\cite[pp.\;241--242]{Je89} (emphasis added)
\end{enumerate}
Ideas such as ``new kind of quantity'' in Leibniz are incompatible
with the Arthur--Rabouin reading.  Jesseph concludes: 
\begin{enumerate}\item[]
In the Leibnizian scheme, true mathematical principles will be found
acceptable on any resolution of the metaphysical problems of the
infinite.  Thus, Leibniz' concern with matters of rigor leads him to
propound a very strong thesis indeed, namely no matter how the symbols
``dx" and ``dy" are interpreted, the basic procedures of the calculus
can be vindicated.  Such vindication could take the form of a new
science of infinity, or it could be carried out along classical lines,
but in either case the new methods will be found completely secure.
\cite[p.\;243]{Je89}
\end{enumerate}
Jesseph's conclusion is consonant with the idea of the presence of two
methods in Leibniz, as argued in the 2013 \emph{Erkenntnis} article
``Leibniz's infinitesimals: Their fictionality, their modern
implementations, and their foes from Berkeley to Russell and beyond''
\cite{13f}.  In fact, Jesseph was roundly criticized by Bassler -- one
of the authors endorsed by Archibald et al.~in their note 3.%
\footnote{See note~\ref{f5}.}
An analysis of Bassler's criticism of Jesseph appears in
\cite[Section~2.3]{21g}.

Jesseph's appeal to two methods in Leibniz re-emerged in his
discussion in \cite{Je15} of Leibniz's method of computing integrals
via transmutation of curves from \emph{De Quadratura Arithmetica},
which requires knowledge of the tangent lines to the curve.  For conic
sections, the tangent lines were known classically, but for the method
to apply more generally, the tangents can only be obtained via ratios
of genuine infinitesimals.  Here, at least, Jesseph endorses genuine
infinitesimals as an irreducible part of the Leibnizian framework; see
further in \cite{16a}.

Similar remarks apply to Archibald's coauthor Panza, who contrasts
Newton's tradition with
\begin{enumerate}\item[]
La deuxi\`eme tradition \ldots\ que j'ai appel\'ee
\emph{infinit\'esimaliste} et qui remonte aux travaux de Leibniz et
Johann I Bernoulli: le \emph{calcul} est consid\'er\'e comme un
algorithme des diff\'erences infiniment petites qui se produisent dans
une certaine quantit\'e lorsqu'une diff\'erence de la m\^eme sorte se
produit dans une quantit\'e li\'ee.%
\footnote{Panza \cite[p.\;xix]{Pa92}; emphasis in the original.}
\end{enumerate}
Since it is generally acknowledged that for Bernoulli, infinitesimals
were mathematical entities (rather than non-`referring' stenography
for exhaustion), Panza's grouping of Leibniz with Bernoulli in his
description of the second tradition puts Panza at odds with the
Ishiguro--Rabouin reading.

To his honor, Craig Fraser (whose earlier critique \cite{Fr15} of
infinitesimal methodology was analyzed in \cite[Sections~4.2--4.6,
  pp.\,123--128]{17a} and \cite[Section~4.4, pp.\;277--278]{18e}) does
not appear among Archibald's coauthors.

\subsection{Non-Archimedean continuum?}

Referring to the argument in the 2013 article \cite{13f}, Archibald et
al.~assert the following:
\begin{enumerate}\item[]
The main claim was that Leibniz shared with Bernoulli a certain view
of the continuum as consisting of infinitesimal numbers in addition to
ordinary (or ``assignable'') numbers.  We may note in passing that
this already involves anachronism at odds with a properly historical
approach.  For Leibniz did not conceive of numbers as constituting a
continuum, nor did he allow infinite sets (infinite wholes, in his
terminology). {\ldots} There is no way that one can claim that
Bernoulli defended a certain picture of the continuum ``following
Leibniz.''  Accordingly, there is no B-methodology sensu stricto in
Leibniz.  \cite[p.\;2]{Ar22b}
\end{enumerate}
Archibald et al.'s objection here is two-fold:
\begin{enumerate}
\item
that Leibniz ``did not conceive of numbers as constituting a
continuum'' and rejected infinite wholes.
\item
that there is no B-methodology (i.e., methodology involving genuine
infinitesimals) in Leibniz.
\end{enumerate}

There are two major problems with Archibald et al.'s claims.

First, using the term \emph{continuum} does not imply either its
punctiform structure or that of an infinite whole, any more so than
does the term \emph{extension}.  Leibniz tended to use the latter
term; he used it for example in the letter to Masson from 1716
analyzed in detail in the article in \emph{Review of Symbolic Logic}
\cite{23d}.  One wonders how Archibald would evaluate the title
\emph{The Labyrinth of the Continuum} of the classical work
\cite{Ar01} edited by {\ldots} his coauthor Arthur.  Archibald et
al.~are therefore attacking a strawman (if not themselves).

Second, Leibniz has two documents from 1695 where he makes it clear
that his incomparables violate the notion of comparability expressed
in Euclid~V definition~4, which is a version of the Archimedean
property; see \cite[p.\;288]{Le95a} and \cite[p.\;322]{Le95b}.  Thus,
Leibniz is rather explicit about non-Archimedean phenomena occurring
for his incomparables.

\subsection{How many tracks?}

Archibald et al.~go on to claim that the 2013 article \cite{13f}
argued that Leibniz had two methods: track A and track B, whereas in
the 2022 ``Two-track depictions,'' Bob asserts on the contrary that
Leibniz used only the track B method whereas Alice claims that Leibniz
used only the track A method:

\begin{enumerate}\item[]
Katz et al. completely changed their position, but without
acknowledging this change, \emph{as if it did not ruin their previous
  argument}.  In the above-cited paper published in this journal,
which is supposed to give a survey of a long-standing debate, A and B
are no longer presented as a pair of methodologies in Leibniz, but as
positions endorsed by commentators to understand the term ``fiction''
in Leibniz.  \cite[p.\;2]{Ar22b} (emphasis added)
\end{enumerate}
Contrary to the claim by Archibald et al., the 2022 article does not
``ruin'' the previous argument at all.  The existence of two methods
in Leibniz is established fact that was commented on in detail by Bos
in the seminal study \cite{Bo74}.  Furthermore, this fact was accepted
by Archibald's coauthor Jesseph; see Section~\ref{jesseph}.  The two
methods are the exhaustion method and the infinitesimal method.  Bos
mentioned that the method using infinitesimals exploited the law of
continuity.%
\footnote{Leibniz's law of continuity possesses an adequate proxy in
  Robinson's transfer principle as explained in \cite{21a}, and has no
  convincing analog in the Weierstrassian setting.}
Furthermore, the existence of two methods in Leibniz is a strong
argument against the Rabouin--Arthur interpretation, which makes it
difficult to distinguish between the two methods, an ``infinitesimal''
being merely stenography for more exhaustion.  The point is that in
\cite{22b}, Alice and Bob are arguing specifically about the
interpretation of Leibniz's term \emph{infinitesimal} (rather than
about Leibniz's exhaustion method which involves no infinitesimals
even nominally).  Thus, we stand by both
\begin{enumerate}
\item
our (and Henk Bos's) position that there are two separate methods in
the Leibnizian calculus, and
\item
our position that Leibniz's infinitesimal method involved genuine
infinitesimals rather than stenography for exhaustion.
\end{enumerate}
There is no contradiction between the two positions, contrary to the
claim by Archibald et al.  The cardinal point here is that Leibniz's
non-infinitesimal (``exhaustion'') method was indeed mentioned by Bob
in the article ``Two-track depictions'' pilloried by Archibald et al.:
\begin{enumerate}\item[]
Bob argues that Archimedean paraphrases in exhaustion style constitute
an alternative method rather than an unwrapping of the infinitesimal
method.  \cite[p.\;262]{22b}
\end{enumerate}

Archibald et al.~are certainly within their rights to disagree with
our arguments, but their attempt to win the argument by
misrepresenting
our position does not amount to a helpful contribution to historical
scholarship.  

In closing, it is ironic that Archibald et al.~should claim that
\begin{enumerate}\item[]
[O]ver the years, it became clearer and clearer that our interlocutors
do not care much about rational discussion and scientific dialogue
from \emph{different perspectives}, but seek rather to disparage their
alleged enemies, {\ldots} The latest example of that approach is
provided by a paper {\ldots}~``Two-Track Depictions of Leibniz’s
Fictions.'' \cite[p.\;2]{Ar22b} (emphasis on ``different
perspectives'' added)
\end{enumerate}
For ``Two-track depictions" is devoted specifically to making explicit
a pair of \emph{different perspectives} on Leibniz's calculus, so as
to stimulate rational discussion and scientific dialogue.

Archibald et al.~do little to clarify the Question of Fundamental
Methodology, namely that the history of mathematics, like mathematics
itself, could benefit from a plurality of approaches.

\section{Conclusion}
\label{s6}

We have examined the modern debate over the infinitesimal calculus
from Leibniz to Cauchy.  Scholars who seek to interpret the Leibnizian
calculus in a purely Archimedean context while denying his
infinitesimal the status of a mathematical entity, often do not
adequately appreciate Leibniz's distinction between infinite wholes
and bounded infinities (the inverse of infinitesimals).  Leibniz's
rejection of the former does not imply rejection of the latter.

To grant Leibniz and Cauchy the use of genuine infinitesimals is not
to impute to them the anticipation of modern nonstandard analysis or
nonstandard numbers.  It only means to argue that, in line with the
procedures/foundations distinction, one finds better proxies for their
ideas and inferential moves in Robinson's framework for analysis with
infinitesimals than in Weierstrassian analysis.

Leibniz and Cauchy had systems of infinitesimal analysis that yielded
correct predictions in analysis, geometry, physics and elsewhere when
one adhered to the internal rules of those systems.

A historiography that wishes to see Leibniz and Cauchy as direct
predecessors of Weierstrass, and therefore denies them any use of
genuine infinitesimals, runs the risk of being teleological.

One should also recognize for mathematics that history is contingent,
in line with the insights by Grattan-Guinness and Hacking, and the
evolution of analysis did not necessarily have to result in the
elimination of infinitesimals.

\section*{Acknowledgments} 

The authors are grateful to Emanuele Bottazzi, Michael Barany, Carlo
Beenakker, Oliver Knill, Karl Kuhlemann, and David Schaps for helpful
comments.  The influence of Hilton Kramer (1928--2012) is obvious.

\end{document}